\documentclass[10pt]{article}
\usepackage{amsmath}
\usepackage{amssymb}
\usepackage{amsfonts}
\usepackage{amsthm}
\usepackage{slashed}
\usepackage{cancel}
\usepackage{mathtools}
\usepackage{paralist}
\usepackage{graphicx}
\usepackage{float}
\usepackage{tikz}
\usepackage{tikz-cd}
\usepackage{pgfplots}
\usepackage[hmargin = 3.7 cm,vmargin = 3.7 cm]{geometry}
\usepackage{hyperref}
\numberwithin{equation}{section}
\usepackage{cleveref}
\usepackage{cite}
\usepackage{array}
\usepackage{setspace}
\usepackage{bbm}
\usepackage{dsfont}
\usepackage{mathrsfs}

\usepackage{ mathrsfs}

\usepackage{blkarray}

\usepgfplotslibrary{fillbetween}
\usetikzlibrary{calc}

\newtheorem{lemma}[equation]{Lemma}
\newtheorem{proposition}[equation]{Proposition}
\newtheorem{theorem}[equation]{Theorem}
\newtheorem{corollary}[equation]{Corollary}
\newtheorem{thm}{Theorem}

\theoremstyle{definition}
\newtheorem{definition}[equation]{Definition}

\theoremstyle{definition}
\newtheorem{notation}[equation]{Notation}

\theoremstyle{remark}
\newtheorem{remark}[equation]{Remark}
\newtheorem*{remark*}{Remark}

\theoremstyle{remark}

\theoremstyle{definition}

\newcommand{\vol}{\mathrm{vol}}
\newcommand{\im}{\mathrm{im}}

\newcommand{\SPAN}{\mathrm{span}}

\renewcommand{\Re}{\operatorname{Re}}
\renewcommand{\Im}{\operatorname{Im}}

\newcommand{\SL}{\mathrm{SL}}

\newcommand{\dlpara}{(\!(}
\newcommand{\drpara}{)\!)}

\newcommand{\dbar}{\bar{\partial}}
\newcommand{\bdry}{|_{M}}
\newcommand{\bdryf}{\|_{\partial,4}}
\newcommand{\bdrytn}{\|_{\partial Z}}

\newcommand{\bC}{{\mathbb{C}}}

\newcommand{\db}{\bar{\partial}}
\newcommand{\bR}{{\mathbb{R}}}
\newcommand{\bZ}{{\mathbb{Z}}}
\newcommand{\bN}{{\mathbb{N}}}

\newcommand{\Lie}{\mathcal{L}}

\newcommand{\tpsi}{\tilde{\psi}}
\newcommand{\cK}{{\cal K}}
\newcommand{\cH}{{\cal H}}
\newcommand{\cV}{{\cal V}}
\newcommand{\cD}{{\cal D}}

\newcommand{\hcK}{\hat{{\cal K}}}

\newcommand{\rcyc}{{\rm cyclic}}

\begin{document}

\title{Calabi-Yau threefolds with boundary}
\author{Simon Donaldson and Fabian Lehmann}
\date{}
\maketitle
\vspace{-16pt}
{\centering{\it In fond memory of Jean-Pierre Demailly}
\vspace{10pt}

\date{\today}

\vspace{12pt}}

In this paper we develop a deformation theory for complex Calabi-Yau threefolds with boundary, parallel to the theory for $G_{2}$-structures developed in \cite{Donaldson-Bdry-G2}. Thus we consider a $C^{
\infty}$ $6$-manifold $Z$ with boundary $M=\partial Z$ and a nowhere-zero holomorphic $3$-form $\Psi+i \hat{\Psi}$, for real forms $\Psi,\hat{\Psi}$. We address the question: if $\Psi\vert_{M}$ is deformed to a nearby closed $3$-form on $M$ can we deform the Calabi-Yau structure on $Z$ to match? In our previous paper \cite{embedding-paper} we studied a variant of this question,  for domains $Z$ in a fixed ambient Calabi-Yau threefold. 

Much of the motivation for this work comes from Hitchin's approach to Calabi-Yau structures in three complex dimensions \cite{Hitchin}, which we recall in outline. On a  $6$-dimensional real vector space $V$ there is a $GL(V)$-invariant open set ${\cal S}\subset \Lambda^{3} V^{*}$ of \lq\lq stable'' elements and a $GL(V)$-equivariant map $P:{\cal S}\rightarrow {\cal S}$ such that each $\Psi\in {\cal S}$ defines a complex structure on $V$ for which $\Psi+ i P(\Psi)$ has type $(3,0)$.  On a $6$-manifold $Z$ a stable $3$-form $\Psi$ (i.e. one which is stable at each point) defines a Calabi-Yau structure if $d\Psi=0, dP(\Psi)=0$, a system of PDE for the real $3$-form
$\Psi$.  Thus the question we address is a boundary value problem for this PDE. 

To set up our problem more precisely we introduce terminology, following \cite{Donaldson-Bdry-G2}.  For our $6$-manifold $Z$ with boundary $M=\partial Z$ we define  an \lq\lq enhanced boundary value'' to be an equivalence class of closed $3$-forms $\Psi$ on $Z$ under the equivalence  $\Psi_{1}\sim \Psi_{2}$ if $\Psi_{1}-\Psi_{2}= d\alpha$ for a $2$-form $\alpha$ on $Z$ with $\alpha\vert_{M}=0$.
Thus an enhanced boundary value defines a closed $3$-form $\psi$  on $M$ and two enhancements of the same $3$-form $\psi$ on $M$ differ by a class in the relative cohomology group  $H^{3}(Z,M)$. 
Of course, such enhancements only exist if the class of $\psi$ in $H^{3}(M)$ extends to $H^{3}(Z)$.

 The more precise version of our problem is: given a closed $3$-form $\psi$ on $M$ and an enhancement can one find a stable $3$-form $\Psi$ in that enhancement class such that $dP(\Psi)=0$?  Thus if $\Psi_{0}$ is any form in the enhancement class this is a PDE for a $2$-form $\alpha$ on $Z$
    $$   d P(\Psi_{0}+d\alpha)=0, $$
with boundary condition $\alpha\vert_{M}=0$ (and the open condition that $\Psi_{0}+d\alpha$ is stable).  

This boundary value problem arises naturally from Hitchin's variational approach. A stable form $\Psi$ defines a volume form ${\rm vol}_{\Psi}$ and the derivative of the volume with respect to variations in 
$\Psi$ is
$$  \delta {\rm vol}_{\Psi} = P(\Psi) \wedge \delta \Psi. $$
We get a volume functional on stable forms in an enhancement class
$$  {\rm Vol}(\Psi)= \int_{Z} {\rm vol}_{\Psi}$$ and our boundary value problem is the Euler-Lagrange equation defining critical points of this functional,  just as in the $G_{2}$ case in \cite{Donaldson-Bdry-G2}. 

Let $\mathcal{G}$ be the group  of diffeomorphisms of $Z$ fixing the boundary pointwise and isotopic to the identity through such diffeomorphisms. Then $\mathcal{G}$ acts on the stable forms in an enhancement class and on the solutions of our boundary value problem.  The model we seek to extend is the case of a closed manifold $Z$, when the analogue of an enhancement class  is just a cohomology class in $H^{3}(Z,\bR)$. In that case, Hitchin showed that a critical point of the volume functional is nondegenerate, modulo diffeomorphisms, from which a version of the local Torelli theorem follows easily: the moduli space of Calabi-Yau structures up to diffeomorphism  is locally identifed with a neighbourhood in $H^{3}(Z;\bR)$. 

Let $\Psi$ be a stable $3$-form defining a Calabi Yau structure on our manifold $Z$ with boundary. The derivative of the map $P$ at $\Psi$ is the linear map $J$ on $3$-forms which, after complexification, acts as $-i$ on $\Omega^{3,0}$ and $\Omega^{2,1}$ and $i$ on $\Omega^{0,2}$ and $\Omega^{1,2}$.  So the composite $dPd$  vanishes on $\Omega^{2,0}$ and $\Omega^{0,2}$ and is equal to $2i\partial\dbar$ on $\Omega^{1,1}$.  We define a vector space $\cK_{\Psi}$ to be the quotient
\begin{equation}   \cK_{\Psi}= \frac{ \{ d\alpha \in \Omega^{3}: \alpha\vert_{M}=0, dJd\alpha=0\}}{ \{ \mathcal{L}_{v}\Psi: v\in \Gamma(TZ), v\vert_{M}= 0\} }. \end{equation}
Here $\mathcal{L}_{v}$ denotes the Lie derivative and the quotient represents the solutions of the linearised equation modulo the linearised diffeomorphism action.  Our main results in this paper  are confined to the case when $Z$ has pseudoconvex boundary (with respect to the complex structure defined by $\Psi$). Then we will show that $\cK_{\Psi}$ is a finite dimensional vector space. We say that $\Psi$ is \lq\lq rigid'' if $\cK_{\Psi}=0$.  Our main result is the following analogue of Hitchin's:
\begin{thm}
\label{MainThm}
Let $\Psi_0$ be the real part of a Calabi-Yau structure on $Z$ with enhanced boundary value $\hat{\psi}_0$.  If $\Psi_0$ is rigid then for any enhanced boundary value  close to $\hat{\psi}_0$    there is a Calabi-Yau structure with real part $\Psi$ close to $\Psi_0$  with that  enhanced boundary value and $\Psi$ is unique in a neighbourhood of $\Psi_0$ up to diffeomorphisms close to the identity which fix the boundary pointwise.  
\end{thm} 
In this statement \lq\lq close'' refers to the $C^{\infty}$ topology for data on the manifold with  boundary and the quotient topology on the set of enhanced boundary values.

We will now outline the content of this paper. In Section \ref{section-background} we recall some  background from complex geometry, analysis and the results of \cite{embedding-paper}. Section \ref{section-gauge} is concerned with setting up the deformation problem. First we need a \lq\lq slice'' to represent the space of all deformations around a solution $\Psi_{0}$ with given enhanced boundary value modulo diffeomorphism. If (for simplicity in this outline) $H^{2}(Z)=0$ we show that this can be taken to be a neighbourhood $U$ of $0$ in
$$    A= \{\alpha\in \Omega^{1,1}: d^{*}\alpha=0, \alpha\bdryf=0\} .$$
Here the notation $\alpha\bdryf=0$ means that on $M$ the components of $\alpha$ in a certain rank $4$ subbundle vanish and  the operator $d^{*}$ is defined using some auxiliary Hermitian metric. The proof that this is a slice is similar to the $G_{2}$ case in \cite{Donaldson-Bdry-G2} but with some significant differences.  These require some variants of the usual Hodge theory on manifolds with boundary.  For any closed stable $3$-form $\Psi$ the $4$-form  $dP(\Psi)$ has type $(2,2)$ with respect to the almost complex structure $I_{\Psi}$.  Thus $dP(\Psi)$ should be viewed as a section of an infinite dimensional vector bundle $\cV$  over the space of closed stable $3$-forms whose fibre $\cV_{\Psi}$ at $\Psi$ is the space of exact forms of type $(2,2)$ with respect to $I_{\Psi}$.  To set up the deformation problem we need to construct a local trivialisation of this bundle. Then our solutions are the zeros of a map  from a neighbourhood $U$ of $0$ in $A$ to the fixed vector space $\cV_{\Psi_{0}}$.  The derivative of this map at $0$ is (assuming, again for simplicity in this outline, that $H^{4}(Z)=0$), 
$$    2 i \partial\dbar : \{\alpha\in \Omega^{1,1}: d^{*}\alpha=0, \alpha\bdryf=0\} \rightarrow \{ \rho \in \Omega^{2,2}: d\rho=0\}. $$
As usual, the crucial thing is to show that if $\Psi_{0}$ is rigid this linear map  is invertible.  This can be formulated in terms of an operator $\cD= 2 i\partial\dbar + *dd^{*}:\Omega^{1,1}\rightarrow \Omega^{2,2}$. One wants to show that for any $\rho \in \Omega^{2,2}$ there is a solution of the linear boundary value problem:
\begin{equation}   \cD\alpha=\rho   \   \   \   \   \  \alpha\bdryf=0\ \ \ \ \ \ \ d^{*} \alpha\vert_{M}=0, 
\label{equ_to_solve}
\end{equation}
with mixed Dirichlet and Neumann boundary conditions.
The differential operator $\cD$ is elliptic. For the analogous linear boundary value problem in the $G_{2}$ case it was shown that the boundary conditions define an elliptic system so standard theory could be applied, leading to a result for the nonlinear problem via the Banach space implicit function theorem. The difference in our current situation is that the boundary conditions in \eqref{equ_to_solve} are not elliptic. This is why we need to restrict to the case of pseudoconvex boundary and we can only hope to find solutions with some loss of derivatives and apply the more sophisticated Nash-Moser theory. 

We have not found in the literature  a general result on hypoelliptic boundary value problems which covers \eqref{equ_to_solve}, and in Section \ref{section-proof} we use a two step argument to produce the inverse. The first step is to solve a 
$\partial\dbar$ equation without the boundary condition $\alpha\bdryf=0$ using the $\dbar$-Neumann theory. The second step is to adjust this solution to satisfy the boundary condition.  In our previous paper \cite{embedding-paper} we introduced a finite dimensional vector space $\cH_{M}$ defined by the restriction of $\Psi_{0}$ to the boundary $M$ of $Z$ and the adjustment in the second step is possible when $\cH_{M}=0$. 

The Nash Moser theory requires more than just the invertibility of the derivative at the given solution. We apply a version of the theory going back to Zehnder \cite{ZehnderIFT} which requires that the derivative is \lq\lq approximately invertible'' at all nearby points. We did not find in the literature an accessible statement of a result which covers our case precisely so in Appendix \ref{appendix-Nash-Moser} we explain the small modifications required to the proof in \cite{StRaymond}.  
At non-integrable stable forms   $\Psi$ close to $\Psi_0$
we cannot invert the derivative and we carefully need to check that the error satisfies the condition from Appendix \ref{appendix-Nash-Moser}. 
We first prove Theorem \ref{MainThm} in the special case when $\mathcal{H}^{1,2}=0$ and $\mathcal{H}_M=0$, and then in the general case $\mathcal{K}_{\Psi_0}=0$. In Section \ref{section-defspace} we make a more detailed analysis of the space $\cK_{\Psi}$ and its relation with
 $\cH_{M}$.  When $Z$ is Stein we find that $\cK_{\Psi}$ is the kernel of a linear map from $\cH_{M}$ to $H^{3}(Z,M)$. We also make a study of the example $Z=B^{3}\times T^{3}$.
 
\vspace{10pt}
This work was supported by the Simons Foundation through the \textit{Simons Collaboration on special holonomy in geometry, analysis and physics}. Much of the work on this paper was carried out while the authors were members of the Simons Center for Geometry and Physics, Stony Brook. The authors are grateful to the Center for its support.
 
 \section{Background}
 \label{section-background}
 
\subsection{Complex Geometry}

Denote by $\Omega^3_s(Z)$ the open subset of stable $3$-forms.
As noted in the introduction there exists a map 
$P:\Omega^3_s(Z)\rightarrow\Omega^3_s(Z)$
such that each $\Psi\in\Omega^3_s(Z)$ induces an almost complex structure $I_{\Psi}$ with respect to which $\Psi+iP(\Psi)$ has type $(3,0)$. The $\SL(3,\bC)$-structure given by $\Psi+iP(\Psi)$ is torsion-free and $I_{\Psi}$ is integrable if and only if
\begin{align}
d\Psi=0, \quad dP(\Psi)=0.
\label{tf-equations}
\end{align}
We have 
\begin{lemma}[see {\cite[Lemma 1]{DonaldsonBdryRemarks}}]
\label{lemma-22}
If $\Psi$ is closed, then $dP(\Psi)$ has type $(2,2)$ with respect to $I_{\Psi}$.
\end{lemma}

We will also need to consider $\Psi$ for which $I_{\Psi}$ is not integrable. 
On a general almost complex manifold the exterior derivative splits as
\begin{align*}
d = \partial + \dbar + i_N.
\end{align*}
If we restrict $d$ to $\Omega^{p,q}$, then $\partial$ is defined as the projection of $d$ to  $\Omega^{p+1,q}$, 
$\dbar$ the projection to  $\Omega^{p,q+1}$, and $i_N$ the anti-derivation of degree $1$ defined by the Nijenhuis tensor $N$ of the almost complex structure. $N\in\Omega^2(Z,TZ)$ splits into two components $N'\in \Omega^{2,0}(Z,T^{0,1}Z)$ and $N''\in\Omega^{0,2}(Z,T^{1,0}Z)$ with $N''=\bar{N'}$.

\begin{lemma}
\label{lemma-N}
If $d\Psi=0$, then $N$ is a function of $\Psi$ and $dP(\Psi)$.
\end{lemma}
\begin{proof}
Suppose that $\theta$ has type $(1,0)$ with respect to $\Psi$, which is equivalent to 
\begin{align*}
\theta\wedge(\Psi+iP(\Psi))=0.
\end{align*}
Applying $d$, decomposing into types and using Lemma \ref{lemma-22}
gives
\begin{align*}
i_{N''}\theta\wedge(\Psi+\sqrt{-1} P(\Psi)) = \sqrt{-1}\, \theta\wedge dP(\Psi).
\end{align*}
\end{proof}

\begin{lemma}
\label{lemma-ACIDs}
$d^2=0$ is equivalent to:
\begin{gather*}
i_{N'}i_{N'}
=0,
\\
\partial i_{N'} + i_{N'}\partial
=0,
\\
\partial^2 + \dbar i_{N'} + i_{N'}\dbar
=0,
\\
\partial\dbar+\dbar\partial +i_{N'}i_{N''}+ i_{N''}i_{N'}
=0,
\\
\dbar^2 + \partial i_{N''}+ i_{N''}\partial
=0,
\\
\dbar i_{N''} + i_{N''}\dbar
=0,
\\
i_{N''}i_{N''}
=0.
\end{gather*}
\end{lemma}

Now let $\Psi$ be a torsion-free $\SL(3,\bC)$-structure  and choose a compatible Hermitian background metric $h$.
We do not require $h$ to be K\"ahler. Thus the metric form $\omega$ might not be closed. Denote by $L$ the Lefschetz-operator which is the wedge product with $\omega$, and by $\Lambda$ its adjoint which is interior multiplication with $\omega$.
The space of real $2$-forms on $Z$ has a decomposition
\begin{align*}
\Omega^2(Z)
=
\Omega^2_1 \oplus \Omega^2_8 \oplus \Omega^2_6
\end{align*}
where the components are given by
\begin{align}
\Omega^2_1 = \Omega^0(Z)\omega,
\quad
\Omega^2_8 = \{\sigma\in \mathrm{Re}\, \Omega^{1,1}: \Lambda\sigma=0  \}
,\quad
\Omega^2_6 = \{i_X\Psi: X\in\Gamma(TZ)\}.
\label{2-forms}
\end{align}
These relate to the decomposition into forms of type $(p,q)$ as
\begin{align*}
\mathrm{Re}\, \Omega^{1,1} = \Omega^2_1\oplus \Omega^2_8,
\quad
\mathrm{Re}\,(\Omega^{2,0}\oplus \Omega^{0,2}) = \Omega^2_6.
\end{align*}
There are also $3$-forms of type $6$:
\begin{align*}
\Omega^3_6 = \{\eta\wedge\omega: \eta\in\Omega^1(Z)\} = L(\Omega^1(Z)).
\end{align*}
Write $d_k, k=1,6,8$, for the projection of $d:\Omega^1(Z)\rightarrow \Omega^2(Z)$ to $\Omega^2_k$
and write $d^3_6$ for the projection of $d:\Omega^2(Z)\rightarrow \Omega^3(Z)$ to $\Omega^3_6$.

\begin{lemma}
\label{lemma-KahlerID}
For $\sigma_k\in\Omega^2_k(Z), k=1,6,8$, we have 
\begin{align*}
d^3_6 \sigma_1 = *Ld^*\sigma_1+a(\sigma),\quad
d^3_6\sigma_6 = \frac{1}{2} *Ld^*\sigma_6+a(\sigma),\quad
d^3_6\sigma_8 = -\frac{1}{2} *Ld^*\sigma_8+a(\sigma),
\end{align*} 
where $a(\sigma)$ stands for a linear operator of order zero in $\sigma$.
\end{lemma}
\begin{proof}
To simplify notation in the course of this proof we will denote any operator which has order zero by $a$.
We will make use of three standard identities.
\begin{enumerate}
\item
$[\Lambda,L]$ acts by multiplication with $3-k$ on $\Lambda^k T^*Z$.
\item
Given $\eta\in\Lambda^1T^*Z$, we have $*L\eta = - LI\eta$. 
\item
The Hermitian identity $[\Lambda,d]=Id^*I + a$, 
where here the operator $I$ acts by multiplication with $i^{p-q}$ on $\Omega^{p,q}$.
\end{enumerate}
If $\mu=L\xi\in\Lambda^3_6 Z$, then
\begin{align*}
\Lambda \mu = \Lambda L\xi = [\Lambda,L]\xi = 2\xi,
\end{align*}
and thus $L\Lambda\mu=2\mu$. This gives for a $2$-form $\sigma$ with $\Lambda\sigma=0$
\begin{align*}
2 d^3_6\sigma = L\Lambda d^3_6\sigma = 
L\Lambda d\sigma=
L[\Lambda,d]\sigma=
 LId^*I\sigma + a(\sigma) = -*Ld^*I\sigma+a(\sigma).
\end{align*} 
$I$ acts as the identity on $\Lambda^2_8$ and minus the identity on $\Lambda^2_6$. The formula follows in these two cases.

Let $\sigma = Lf$. Then 
\begin{align*}
[\Lambda,d]\sigma
=
\Lambda L df - d \Lambda L f
+a(\sigma)
=
[\Lambda,L]df-d [\Lambda,L]f
+a(\sigma)
=
2df-3df
+a(\sigma)=-df+a(\sigma),
\end{align*}
and thus 
\begin{align*}
d\sigma = Ldf + a(\sigma) = -L[\Lambda,d]\sigma +a(\sigma) = -L I(d^*\sigma)+a(\sigma) = * L d^*\sigma+a(\sigma).
\end{align*}
\end{proof}

\begin{lemma}
\label{lemma-d*d-formula}
For $\eta\in\Omega^1(Z)$ we have
\begin{align*}
d^*d\eta = -d^*d_1\eta + 2d^*d_8\eta+ D\eta,
\end{align*}
where $D$ denotes a linear differential operator of order one.
\end{lemma}
\begin{proof}
In this proof we denote any operator of order one by $D$.
Applying Lemma \ref{lemma-KahlerID} with $\sigma_k=d_k\eta$ gives 
\begin{align*}
0 = d^3_6d\eta
=
d^3_6d_1\eta + d^3_6d_8\eta+d^3_6d_6\eta
=
*L d^*d_1\eta -\frac{1}{2}*Ld^*d_8\eta+\frac{1}{2}*Ld^*d_6\eta
+D\eta.
\end{align*}
Applying the inverse of the isomorphism $*\circ L: \Omega^1(Z)\rightarrow \Omega^3_6$, we get
\begin{align*}
d^*d_6\eta = -2 d^*d_1\eta+ d^*d_8\eta+D\eta.
\end{align*}
Subsituting this into
\begin{align*}
d^*d\eta = d^*d_1\eta + d^*d_8\eta+d^*d_6\eta
\end{align*}
gives the result.
\end{proof}

\subsection{Some elliptic boundary value problems}

By choosing an atlas for $Z$ and a subordinate partition of unity
we can define as usual Sobolev spaces $L^2_s$ for sections of $TZ$ and its associated bundles. We denote the norm of $L^2_s$ by $\|\cdot \|_s$. 
Here $s$ is the number of derivatives if it is an integer, but we also need to consider non-integral $s$. For details on fractional Sobolev spaces we refer to \cite[Appendix 1. and 2.]{FollandKohn}.
We also use $\mathcal{C}^k$-norms which we denote by $[[\,\cdot\, ]]_k$. Here $k$ is an integer.

We first recall the classical Hodge decompositions on Riemannian manifolds with boundary. We need the following spaces:
\begin{align}
\mathcal{H}^k 
&=
\{
\gamma\in \Omega^k(Z): d\gamma = 0, d^*\gamma=0
\},
\nonumber
\\
\mathcal{H}_D^k 
&=
\{
\gamma\in \Omega^k(Z): d\gamma = 0, d^*\gamma=0, \gamma\bdry=0
\},
\nonumber
\\
\mathcal{E}^k_D 
&=
\{d\alpha: \alpha\in\Omega^{k-1}(Z), \alpha\bdry = 0 \},
\label{exact-forms}
\\
\mathcal{C}^k_D 
&= \{ d^*\beta: \beta\in\Omega^{k+1}(Z), \beta\bdry=0\},
\nonumber 
\\
\mathcal{C}^k_N 
&= \{ d^*\beta: \beta\in\Omega^{k+1}(Z), *\beta\bdry=0\}.
\nonumber
\end{align}
Here the subscripts $D$ and $N$ indicate the Dirichlet and von Neumann boundary conditions, respectively.

\begin{proposition}
\cite[Theorem 2.4.2]{Schwarz-ellBVP}
\label{Hodge-mixed}
There is an $L^2$-orthogonal decomposition
\begin{align*}
\Omega^k(Z)
=
\mathcal{H}^k
\oplus
\mathcal{E}_D^k
\oplus
\mathcal{C}_N^k.
\end{align*}
Consequences of this decomposition are:
\begin{compactenum}[(i)]
\item
The space of closed $k$-forms is equal to $\mathcal{H}^k
\oplus
\mathcal{E}_D^k$.
\item 
$\mathcal{C}^k_N$ is the $L^2$-orthogonal complement of the space of closed $k$-forms.
\end{compactenum}
\end{proposition}

\begin{proposition}
\cite[Theorem 2.2.6]{Schwarz-ellBVP}
\label{HodgeDirichlet}
The space $\mathcal{H}_D^k$ is finite dimensional and  there is an $L^2$-orthogonal decomposition
\begin{align*}
\Omega^k(Z)
=
\mathcal{H}_D^k
\oplus
\mathcal{E}_D^k
\oplus
\mathcal{C}_D^k.
\end{align*}
Consequences of this decomposition are:
\begin{compactenum}[(i)]
\item
The space of co-closed $k$-forms is equal to $\mathcal{H}^k_D
\oplus
\mathcal{C}_D^k$.
\item
$\gamma\in\Omega^k(Z)$ is co-exact if and only if 
$\gamma\perp \mathcal{H}_D^k
\oplus
\mathcal{E}_D^k$.
\end{compactenum}
\end{proposition}

\begin{lemma}
\label{lemma-EllBdry}
Let $\lambda\in \Omega^1(Z)$. 
Consider the boundary value problem to find $\eta\in\Omega^1(Z)$
which satisfies 
\begin{subequations}
\label{ellBVP11}
\begin{align}
(d^*d_{1,1} + dd^*)\eta &= \lambda,
\\
\eta|_{\partial Z} &=0,
\label{ellBVP11-BC1}
\\
(d^*\eta)|_{\partial Z}&=0.
\label{ellBVP11-BC2}
\end{align}
\end{subequations}
Let $\mathcal{K}$ be the space of solutions for $\lambda=0$.
\begin{compactenum}[(i)]
\item $\mathcal{K}$ is finite dimensional and \eqref{ellBVP11} has a solution if and only if $\lambda$ is $L^2$-orthogonal to $\mathcal{K}$.
\item If $\eta\in\mathcal{K}$, then $d_{1,1}\eta=0$ and $d^*\eta=0$.
\item If $\lambda\in d^*\Omega^{1,1}$ then there exists a solution, and every solution is co-closed.
\end{compactenum}
\end{lemma}
\begin{proof}
Suppose that $\eta\in\Omega^1(Z)$ satisfies the boundary conditions in \eqref{ellBVP11}. Then we can integrate by parts and use Lemma \ref{lemma-d*d-formula} to obtain
\begin{gather*}
(d^*d_{1,1}\eta,\eta) + (dd^*\eta,\eta)
=
\|d_1\eta\|^2 + \|d_8\eta\|^2 + \|d^*\eta\|^2
\geq -\frac{1}{2} \|d_1\eta\|^2 + \|d_8\eta\|^2 + \|d^*\eta\|^2
\\
=
\frac{1}{2}(-d^*d_1\eta +2  d^*d_8\eta,\eta) + \|d^*\eta\|^2
=
\frac{1}{2}(d^*d\eta,\eta)+\|d^*\eta\|^2 -\frac{1}{2}(D\eta,\eta)
\\
\gtrsim \|\eta\|^2_1-(sc)\|\eta\|_1^2 -(lc)\|\eta\|,
\end{gather*}
where in the last line we have used the coercive estimate for the standard Laplace operator with Dirichlet boundary data.
Rearranging yields the coercive estimate
\begin{align*}
\|\eta\|_1^2 \lesssim
\|d_{1,1}\eta\|^2 + \|d^*\eta\|^2 + \|\eta\|^2.
\end{align*}
(i) follows from this coercive estimate 
as in \cite[Sections 2.2-2.4]{Schwarz-ellBVP} for the standard Laplacian by defining a potential and proving its regularity.
\\
(ii) follows from the  first line in the above estimate.
To prove (iii), let $\gamma\in\Omega^{1,1}$ and $\eta\in\mathcal{K}$.
Then because of \eqref{ellBVP11-BC1} we have
\begin{align*}
(d^*\gamma,\eta) = (\gamma , d\eta) = (\gamma, d_{1,1}\eta).
\end{align*}
This expression vanishes by (ii) and existence of a solution follows from (i).
If we have a solution
\begin{align*}
d^*d_{1,1}\eta + dd^*\eta = d^*\gamma,
\end{align*}
then $dd^*\eta=0$, as this term by \eqref{ellBVP11-BC2} is $L^2$-orthogonal to the other two terms. Integrating by parts again gives $d^*\eta=0$.
\end{proof}

\begin{lemma}
\label{lemma-6BVP}
For every $\kappa\in\Omega^2_6(Z)$, there is a unique $\zeta=G_6(\kappa)\in\Omega^2_6(Z)$ which solves the  boundary value problem
\begin{align*}
\pi_6 d^*d\zeta = \kappa,
\\
\zeta\bdry=0.
\end{align*}
Furthermore, for each $n$ we have the elliptic estimate
\begin{align}
\|G_6(\kappa)\|_{n+2}
\lesssim
\|\kappa\|_n.
\label{constant6}
\end{align}
\end{lemma}
\begin{proof}
The boundary condition allows to integrate by parts, so the problem is self-adjoint. 
Given $\zeta\in\Omega^2_6$ with $\zeta\bdry=0$, the coercive estimate for the standard Laplacian and Lemma \ref{lemma-KahlerID} give the coercive estimate
\begin{align*}
\|\zeta\|_1^2 
&\lesssim
\|d\zeta\|^2 + \|d^*\zeta\|^2 + \|\zeta\|^2
\lesssim 
\|d\zeta\|^2 + \|d^3_6\zeta\|^2 + \|\zeta\|^2
\\
&\lesssim
\|d\zeta\|^2  + \|\zeta\|^2
= (\pi_6 d^*d\zeta,\zeta) + \|\zeta\|^2.
\end{align*}
Suppose that $\zeta = i_w \Psi$ is a solution for $\kappa=0$. Then 
$\mathcal{L}_w\Psi =0$ and $w\bdry=0$.
Then $w$ is in particular holomorphic, and because it vanishes on the boundary it has to vanish everywhere. 
Thus for $\kappa=0$ there is only the trivial solution.
For general $\kappa$ existence and uniqueness of a solution follows again like for the standard Laplacian with Dirichlet boundary data \cite[Sections 2.2-2.4]{Schwarz-ellBVP}.
\end{proof}

\subsection{Review of the $\dbar$-Neumann problem}

\label{section-dbar}

In this subsection let $Z$ be a compact, complex manifold of complex dimension $n$ with boundary $\partial Z$. Write $I$ for the complex structure and let $h$ be a Hermitian metric. Denote by $\dbar^*$ the formal adjoint of $\dbar$ with respect to $h$,
and by $\Delta_{\dbar}$ the $\dbar$-Laplacian.
Given $\gamma\in\Omega^{p,q}(Z)$, the $\dbar$-Neumann problem is to find $\phi\in\Omega^{p,q}(Z)$ which solves the boundary value problem
\begin{subequations}
\label{dbar-Neumann-prob}
\begin{gather}
\Delta_{\dbar}\phi = \gamma, 
\\
p(\phi) = 0,
\label{bdry1} 
\\
q(\phi)=0,
\label{bdry2}
\end{gather}
\end{subequations}
where the boundary operators are given by 
$p = \sigma(\dbar^*,\nu)$ and $q = \sigma(\dbar^*,\nu)\circ \dbar$.
Here $\sigma(\dbar^*,\nu)$ denotes the principal symbol of $\dbar^*$ in the direction of the unit co-normal $\nu$ at the boundary. If $N$ is a unit normal, then $\sigma(\dbar^*,\nu)$ is the interior multiplication with the $(0,1)$-part of $N$. 
These boundary conditions are chosen such that we can integrate by parts:
\eqref{bdry1} is equivalent to 
\begin{align*}
(\dbar^*\phi,\psi) = (\phi, \dbar \psi)
\end{align*}
for all $\psi\in\Omega^{p,q-1}$,
and \eqref{bdry2} is equivalent to
\begin{align*}
(\dbar^*\dbar \phi, \gamma) = (\dbar\phi,\dbar\gamma)
\end{align*}
for all $\gamma\in\Omega^{p,q}$.

The analytic properties of the $\dbar$-Neumann problem rely on the geometry of the boundary. 
$Z$ is called \textit{strongly pseudoconvex}, if at each point of the boundary the Levi form is positive definite. 
On a strongly pseudoconvex complex manifold, 
the $\dbar$-Neumann problem is sub-elliptic. This means that
there is the estimate
\begin{align*}
\int_{\partial Z} |\phi|^2 \vol_{\partial Z}
\lesssim \|\dbar\phi\|^2 + \|\dbar^*\phi\|^2 + \|\phi\|^2
\end{align*}
for all $\phi$ which satisfy the boundary conditions \eqref{bdry1} and \eqref{bdry2}.
\begin{proposition}
\cite{FollandKohn}
\label{prop-dbar-Neumann}
Suppose that $Z$ has a strongly pseudoconvex boundary.
Then we have: 
\begin{itemize}
\item
The operator $\Delta_{\dbar}$ is hypo-elliptic, i.e. if $\phi$ is a distributional solution of the equation $\Delta_{\dbar}\phi=\gamma$ and $\gamma$ is smooth, then $\phi$ is smooth.
\item The space
\begin{align*}
\mathcal{H}^{p,q} = \{\phi\in\Omega^{p,q}(Z): \dbar\phi=0, \dbar^*\phi=0, p(\phi)=0 \}
\end{align*}
is finite dimensional.
\item \eqref{dbar-Neumann-prob} has a solution if and only if $\gamma$  is $L^2$-orthogonal to $\mathcal{H}^{p,q}$. In this case we denote by $G_{\dbar}\gamma$ the unique solution orthogonal to $\mathcal{H}^{p,q}$. 
\item If $\dbar\gamma=0$ and $\gamma\perp\mathcal{H}^{p,q}$, then $\beta= \dbar^* G_{\dbar}\gamma$ is a solution of the equation
\begin{align*}
\dbar \beta = \gamma.
\end{align*}
\end{itemize}
\end{proposition}

On a Calabi--Yau $n$ fold with boundary, we have a vanishing result for $\dbar$-harmonic forms of type $(0,n)$.

\begin{proposition}
\label{prop-top-degree-Neumann}
Suppose $Z$ has a holomorphic non-vanishing $(n,0)$-form. Then $\mathcal{H}^{0,n}=0$. 
\end{proposition}
\begin{proof}
For $\gamma\in\Omega^{0,n}(Z)$, the boundary condition $\sigma(\dbar^*,\nu)\gamma=0$ is equivalent to $\gamma$ vanishing as a section of $\Lambda^{0,n}Z$ on the boundary. Thus on $\Omega^{0,n}(Z)$ the $\dbar$-Neumann boundary condition  is elliptic without any assumption on the Levi form. If $Z$ has a Calabi--Yau structure $\Psi+i\hat{\Psi}$, then for $\gamma= f(\Psi-i\hat{\Psi})$ the condition $\dbar^*\gamma=0$ is equivalent to $\bar{f}$ being a holomorphic function. Because $f$ has to vanish on the boundary, by the maximum principle we get $f\equiv 0$.
\end{proof}

Next we turn to the $\dbar$-Neumann problem with inhomogeneous boundary data. For points $z\in \partial Z$ define
\begin{align}
\Lambda^{p,q}_T|_z
:=
\{
\phi\in \Lambda^{p,q} Z|_z : p(\phi)=0
\}.
\label{LambdaT}
\end{align}

\begin{lemma}
\label{lemma-inhomo-bdry-data}
For every $\gamma\in\Omega^{p,q}$, $\xi\in\Gamma(\partial Z, \Lambda^{p,q-1}_T)$ and $\rho \in\Gamma(\partial Z, \Lambda^{p,q}_T)$, there exist unique $\phi\perp \mathcal{H}^{p,q}$ and $h\in\mathcal{H}^{p,q}$ which solve
\begin{align*}
\Delta_{\dbar}\phi + h = \gamma, 
\quad
p(\phi) = \xi,
\quad
q(\phi)=\rho.
\end{align*}
\end{lemma}
\begin{proof}
Uniqueness is clear. To prove existence, we first show that we can find $\varphi\in\Omega^{p,q}(Z)$ such that
\begin{subequations}
\begin{gather}
p(\varphi) = \xi,
\label{solve-bdry-1}
\\
q(\varphi) = \rho.
\label{solve-bdry-2}
\end{gather} 
\end{subequations}
Choose a collar neighbourhood of $\partial Z$ in $Z$. We then have a radial function $r$ such that $\partial_r|_{\partial Z} = N$, where $N$ is a inward pointing unit normal. Write
\begin{align*}
\varphi
=
\dbar r\wedge \alpha + \beta,
\end{align*}
where $\partial_r^{0,1}\lrcorner\beta=0$.
Then \eqref{solve-bdry-1} is equivalent to $\alpha|_{\partial Z} = \xi$.
Writing $L = \dbar - \dbar r \wedge \partial_r^{0,1}$, equation \eqref{solve-bdry-2} is equivalent to
\begin{align}
N\cdot\beta = i IN\cdot \beta + 2\rho + 2L\alpha + A(\alpha,\beta),
\label{solve-bdry-2-explicit}
\end{align}
where $A(\alpha,\beta)$ denote a zero order term in $\alpha$ and $\beta$. $L$ only involves tangential derivatives on the boundary, and $IN$ is tangential as well. Thus any choice of $\beta|_{\partial Z}$ determines a right hand side in \eqref{solve-bdry-2-explicit}, and then $\beta$ can be extended to the interior to solve \eqref{solve-bdry-2}.

We can assume that $\varphi$ is orthogonal to $\mathcal{H}^{p,q}$.
By Proposition \ref{prop-dbar-Neumann} there are  $\phi' \perp \mathcal{H}^{p,q}$ and $h\in\mathcal{H}^{p,q}$
which solve
 \begin{gather*}
\Delta_{\dbar}\phi' + h = \gamma-\Delta_{\dbar}\varphi, 
\quad
p(\phi') = 0,
\quad
q(\phi')=0.
\end{gather*}
Then $\phi := \varphi + \phi'$ is the desired solution.
\end{proof}

To apply Nash--Moser theory, we have to solve a $\dbar$-Neumann problem for the $\dbar$-operator associated with all the different almost complex structures induced by the $\SL(3,\bC)$-structures $\Psi$ in a neighbourhood of a reference structure $\Psi_0$.
Almost complex structures close to a reference structure $I_0$ are parametrised by linear bundle maps $\mu:T^{0,1}\rightarrow T^{1,0}$ by setting 
$T^{0,1}_{\mu}=\mathrm{graph}(-\mu)=\{w-\mu w: w\in T^{0,1} \}$. 
The map $\mu$ has a conjugate
$\bar{\mu}:T^{1,0}\rightarrow T^{0,1}$. Then $T^{1,0}_{\mu}=\mathrm{graph}(-\bar{\mu})$. We also have dual maps 
$\mu^*:\Lambda^{1,0}\rightarrow \Lambda^{0,1}$ and $\bar{\mu}^*:\Lambda^{0,1}\rightarrow \Lambda^{1,0}$. Then $\Lambda^{1,0}_{\mu}=\mathrm{graph}(\mu^*)$ and $\Lambda^{0,1}_{\mu}=\mathrm{graph}(\bar{\mu}^*)$.

The inverse of the the isomorphism $\mathbbm{1}-\mu:T^{0,1}\rightarrow T^{0,1}_{\mu}$ is the projection $\pi^{0,1}:T^{0,1}_{\mu}\rightarrow T^{0,1}$. The dual of this map gives an isomorphism $c_{\mu}:\Lambda^{0,1}\rightarrow \Lambda^{0,1}_{\mu}$.
Taking conjugates and wedge products of $c_{\mu}$ gives
bundle isomorphisms
\begin{align*}
c_{\mu}: \Lambda^{p,q}\rightarrow \Lambda^{p,q}_{\mu}.
\end{align*}
Write $I_{\mu}$ for the almost complex structure determined by $\mu$ and $\dbar_{[\mu]}$ for the $\dbar$-operator associated with $I_{\mu}$. We allow any $\mu$. In particular, $I_{\mu}$ does not need to be integrable.

Let $h_{\mu}$ be a family of metrics which varies smoothly in $\mu$
such that $h_{\mu}$ is Hermitian for $I_{\mu}$. Write $\dbar_{[\mu]}^*$ for the formal adjoint of $\dbar_{[\mu]}$ with respect to $h_{\mu}$ and set $\Delta_{\dbar}(\mu):= \dbar_{[\mu]}\dbar_{[\mu]}^*+\dbar_{[\mu]}^*\dbar_{[\mu]}$.
We want to study the $\dbar$-Neumann problem with respect to $I_{\mu}$.
We prefer to study operators between fixed bundles.
Define 
\begin{align*}
\dbar_{\mu}:= c_{\mu}^{-1} \circ \dbar_{[\mu]}\circ c_{\mu}
\end{align*}
 to be the pull-back of the operator $\dbar_{[\mu]}:\Omega^{p,q}_{\mu}\rightarrow \Omega^{p,q+1}_{\mu}$ to obtain an operator $\Omega^{p,q}\rightarrow \Omega^{p,q+1}$. In other words, we have a commutative diagram
\begin{equation}
\begin{tikzcd}
\Omega^{p,q} \ar[r,"\dbar_{\mu}"] \ar[d,"c_{\mu}"'] & \Omega^{p,q+1} \ar[d,"c_{\mu}"]
\\
\Omega^{p,q}_{\mu} \ar[r,"\dbar_{[\mu]}"] & \Omega^{p,q+1}_{\mu}
\end{tikzcd}
\end{equation}

In the following we want to derive a formula for $\dbar_{\mu}$. 
 Suppose that $(z^1,\cdots, z^n)$ is a holomorphic chart for the complex structure $I_0$. 
For an increasing multi-index $A=(\alpha_1, \cdots ,\alpha_q)$, write 
$d\bar{z}^A$ for $dz^{\alpha_1}\wedge \cdots \wedge d\bar{z}^{\alpha_p}$. Using the summation convention, in this coordinate system we can write
\begin{align*}
\varphi=\varphi_{A,I}\, d\bar{z}^A \otimes dz^I,
\quad
\mu = \mu^{\alpha}_{\beta}\, d\bar{z}^{\beta}\otimes \frac{\partial}{\partial z^{\alpha}.}
\end{align*}
Define $\varepsilon^A_B$ to be $\pm 1$ if $A$ is a permutation of $B$ of sign $\pm 1$, and define $\varepsilon^A_B$ to be $0$ otherwise.
Define a first order differential operator $D_{\mu}: \Omega^{p,q}\rightarrow \Omega^{p,q+1}$ by
\begin{align*}
D_{\mu}\varphi_{C,I}=\varepsilon_C^{\theta A} \mu^{\rho}_{\theta} 
\frac{\partial \varphi_{A,I}}{\partial z^{\rho}}.
\end{align*}
A calculation shows that $D_{\mu}$ is coordinate invariant.

The next Lemma is analogous to Theorem 1 in \cite{Hamilton1977Deformation1}.
\begin{lemma}
For every $\mu$ there exists a tame operator of order zero $a_{\mu}$  such that  we have
\begin{align}
\dbar_{\mu}
=
\dbar - D_{\mu} + a_{\mu}
.
\label{dbar-relation}
\end{align}
\end{lemma}
\begin{proof}
Let $p\in Z$ be arbitrary.
Take coordinates $(z^1,\cdots, z^n)$ centered in $p$ which are holomorphic for $I_0$ and $(W_1,\cdots, W_n)$ a frame for $T^{1,0}_{\mu}$ in a neighbourhood of $p$ with dual frame $\omega^1,\cdots, \omega^n$. We can choose $n$ complex functions  $w^1, \cdots, w^n$ which form a coordinate system such that at $p$ we have $\frac{\partial}{\partial w^{\alpha}} = W_{\alpha}$, and in particular that $T^{0,1}_{\mu}$ annihilates the functions $w^{\beta}$ in $p$. Then the calculation at the beginning of the proof of Theorem 1 in \cite{Hamilton1977Deformation1} is still valid at the point $p$, so that at $p$ we have
\begin{align}
\overline{W}_{\gamma}
=
\frac{\partial \bar{z}^{\beta}}{\partial \bar{w}^{\gamma}}
\left(
\frac{\partial}{\partial \bar{z}^{\beta}}-\mu^{\eta}_{\beta} \frac{\partial}{\partial z^{\eta}}
\right).
\label{fundamentalequ}
\end{align}
Introduce the transition functions $T^{\alpha}_{\beta}$ defined by
\begin{align*}
c_{\mu}(d\bar{z}^{\alpha}) = T^{\alpha}_{\beta} \bar{\omega}^{\beta}. 
\end{align*}
Writing 
\begin{align*}
T^A_B
=
\sum_{
\pi\in S_q}
(\text{sgn}\, \pi) T^{\pi(\alpha_1)}_{\beta_1} \cdots T^{\pi(\alpha_q)}_{\beta_q}
\end{align*}
we then have for $\psi=c_{\mu}\varphi$
\begin{align*}
\psi_{B,J}= T^A_B \overline{T}^I_J \varphi_{A,I}.
\end{align*}
From \eqref{fundamentalequ} it follows that
\begin{align*}
T^{\alpha}_{\beta}(p) = \frac{\partial \bar{z}^{\alpha}}{\partial \bar{w}^{\beta}}(p),
\end{align*}
so that at $p$ we really have
\begin{align*}
\overline{W}_{\gamma}
=
T^{\beta}_{\gamma}
\left(
\frac{\partial}{\partial \bar{z}^{\beta}}-\mu^{\eta}_{\beta} \frac{\partial}{\partial z^{\eta}}
\right).
\end{align*}
In the following, denote by l.o.t. lower order terms. We have
\begin{align*}
\dbar_{[\mu]}\psi_{C,J}
&=
\varepsilon^{\gamma B}_{C} \overline{W}_{\gamma}\psi_{B,J}
+
\text{l.o.t.}
\\
&=
\varepsilon^{\gamma B}_{C} T^A_B \overline{T}^I_J
\overline{W}_{\gamma}\varphi_{A,I}+\text{l.o.t.}
\\
&=
\varepsilon^{\gamma B}_{C} T^A_B \overline{T}^I_J
T^{\theta}_{\gamma}
\left(\frac{\partial\varphi_{A,I}}{\partial \bar{z}^{\theta}}
-
\mu^{\eta}_{\theta} \frac{\partial\varphi_{A,I}}{\partial z^{\eta}}
\right)
+\text{l.o.t.}
\end{align*}
Now we have
\begin{align*}
\varepsilon^{\gamma B}_{C} T^A_B 
T^{\theta}_{\gamma}
=
\varepsilon^{\theta A}_{D} T^D_C,
\end{align*}
which gives 
\begin{align*}
(c_{\mu}\dbar_{\mu}\varphi)_{C,J}
&=
\dbar_{[\mu]}\psi_{C,J}
\\
&=
 T^D_C \overline{T}^I_J
\left(\varepsilon^{\theta A}_{D}\frac{\partial\varphi_{A,I}}{\partial \bar{z}^{\theta}}
-
\varepsilon^{\theta A}_{D}
\mu^{\eta}_{\theta} \frac{\partial\varphi_{A,I}}{\partial z^{\eta}}
\right)
+\text{l.o.t.}
\\
&=
(c_{\mu}(\dbar-D_{\mu})\varphi)_{C,J}  + \text{l.o.t.}
\end{align*}
\end{proof}
\begin{remark*}
A more refined calculation like in the proof of Theorem 1 in \cite{Hamilton1977Deformation1} shows that $a_{\mu}$ vanishes if $I_{\mu}$ is integrable, but this will not be important for our estimate.
\end{remark*}

If $L_1, \cdots, L_n$ is a local frame for $T^{1,0}$, then set
\begin{align*}
\bar{L}_c^{\mu} = \bar{L}_c - \mu^a_c L_a.
\end{align*}
Then we have
\begin{align}
\dbar_{\mu}\varphi_{A,I} = \varepsilon^{cC}_A \bar{L}^{\mu}_c\varphi_{C,I}+\text{l.o.t.}
\label{loc-formula-mu-op}
\end{align}
Denote by $(\cdot, \cdot)_{\mu}$ the  $L^2$-product associated with $h_{\mu}$. 
We can use the isomorphisms $c_{\mu}$ to pull back $h_{\mu}$ to obtain a smooth family of Hermitian metrics on $\Lambda^{p,q}$ and we denote the associated $L^2$-product by $\dlpara\cdot,\cdot\drpara_{\mu}$.
The boundary condition \eqref{bdry1} for $\psi\in\Omega^{p,q}_{\mu}$ is equivalent to
\begin{align*}
(\dbar^*_{[\mu]}\psi,\phi)_{\mu} = (\psi,\dbar_{[\mu]}\phi)_{\mu}
\end{align*}
for all $\phi\in\Omega^{p,q-1}_{\mu}$. 
Thus under the isomorphism $c_{\mu}$ this boundary condition on $\varphi\in\Omega^{p,q}$ is
\begin{align*}
\dlpara \dbar^*_{\mu}\varphi,\phi\drpara_{\mu} = \dlpara\varphi, \dbar_{\mu}\phi\drpara_{\mu}
\end{align*} 
for all $\phi\in\Omega^{p,q-1}$. 

Hamilton considers the $\dbar$-complex for $(0,q)$-forms with values in $T^{1,0}$ while we consider the $\dbar$-complex for $(0,q)$-forms with values in $\Lambda^{p,0}$. However, at leading order the formula \eqref{loc-formula-mu-op} is analogous to the formula for Hamilton's $\dbar_{\mu}$-operator \cite[p.18]{Hamilton1977Deformation1}. Thus all his results carry over in a straightforward way. In particular, Hamilton \cite[p.19]{Hamilton1977Deformation1} shows that the Hermitian metrics $h_{\mu}$ can be chosen in such a way that the first $\dbar$-Neumann boundary condition is independent of $\mu$, i.e. 
\begin{align*}
N^{0,1}\lrcorner \varphi = 0
\end{align*} 
where $N$ is a unit normal for $h_0$.
Then the proof of the sub-elliptic estimate in \cite{Hamilton1977Deformation1} carries over without any changes.

\begin{proposition}
\label{prop-uniformMorrey}
Suppose that $Z$ is strongly pseudoconvex with respect to $I_0$. Then for
all $\mu$ in a neighbourhood of $0$
and  all $\varphi\in\Omega^{p,q}(Z)$  
with $p(\varphi)=0$ and $q_{\mu}(\varphi)=0$
we have 
\begin{align*}
\int_{\partial Z} |\varphi|^2 \vol_{\partial Z}
\lesssim \|\dbar_{\mu}\varphi\|_{\mu}^2
+
\|\dbar^*_{\mu}\varphi\|_{\mu} + \|\varphi\|^2
.
\end{align*}
\end{proposition}
Set $E_{\mu} = c_{\mu}^{-1}\circ \Delta_{\dbar}[\mu]\circ c_{\mu}$,
and define a bilinear form
\begin{align*}
Q_{\mu}(\varphi,\psi) = \dlpara E_{\mu}\varphi,\psi\drpara_{\mu}.
\end{align*}
Then by Proposition \ref{prop-uniformMorrey}, $Q_{\mu}$
satisfies  what Hamilton \cite[p.437]{Hamilton1977Deformation2} calls the uniform persuasive estimate
\begin{align*}
\int_{\partial Z} |\varphi|^2 \vol_{\partial Z}
\lesssim \text{Re}\, Q_{\mu}(\varphi,\varphi) + \|\varphi\|^2
\end{align*}
if $\varphi$ satisfies the corresponding boundary conditions.
Thus we can apply the theory developed in \cite[Part 4]{Hamilton1977Deformation2}.

\begin{proposition}
\cite[p. 452]{Hamilton1977Deformation2}
\label{Prop-Apriori}
There exists $s\in\bN$ such that for all $\mu$ in a neighbourhood of zero we have the uniform a priori estimate 
\begin{align}
\|\varphi\|_n
\lesssim
\|E_{\mu}\varphi\|_n
+
|p(\varphi)|_{n+2}
+
|q_{\mu}(\varphi)|_{n+1}
+
(\|\mu\|_{n+s}+1)\|\varphi\|.
\label{apriori}
\end{align} 
\end{proposition}

\begin{lemma}
\cite[Lemma on p.453]{Hamilton1977Deformation2}
\label{lemma-perpToKer}
There exists $l$ such that for all $\mu$ in a neighborhood of $0$ and  all $\varphi \perp \mathcal{H}^{p,q}$ we have 
\begin{align*}
\|\varphi\| \lesssim \|E_{\mu}\varphi\|_l + |p(\varphi)|_{l+2} + |q_{\mu}(\varphi)|_{l+1}.
\end{align*}
\end{lemma}

\begin{proposition}
\label{Prop-UniSolveH0}
For each $\mu$ in some neighborhood of $0$ and for each $\gamma\in\Omega^{p,q}$, $\xi\in\Gamma(\partial Z,\Lambda^{p,q-1}_T)$
and $\rho \in\Gamma(\partial Z,\Lambda^{p,q}_T)$, where $\Lambda^{p,q}_T$ is the bundle defined in \eqref{LambdaT}, there are unique $\phi\perp \mathcal{H}^{p,q}$ and $x\in\mathcal{H}^{p,q}$ such that
\begin{align*}
E_{\mu}\phi + x =\gamma,
\quad
p(\phi) = \xi,
\quad
q_{\mu}(\phi) = \rho.
\end{align*}
There exist $b$ and $s$ such that $(\phi,x)$ satisfies for any $n$ the tame estimate
\begin{align*}
\|\phi\|_n + |x|
\lesssim
\|\gamma\|_n
+
|\xi|_{n+2}
+
|\rho|_{n+1}
+
(\|\mu\|_{n+s}+1)
(\|\gamma\|_b
+
|\xi|_{b+2}
+
|\rho|_{b+1}).
\end{align*}
\end{proposition}
\begin{proof}
Consider the $\mu$-dependent map
\begin{gather*}
(\mathcal{H}^{p,q})^{\perp} \oplus \mathcal{H}^{p,q}
\rightarrow 
\Omega^{p,q} \oplus \Gamma(\partial Z,\Lambda^{p,q-1}_T) \oplus \Gamma(\partial Z,\Lambda^{p,q}_T),
\\
(\phi,x)
\mapsto
(E_{\mu}\phi + x, p(\phi), q_{\mu}(\phi)).
\end{gather*}
For $\mu=0$ this is invertible by Lemma \ref{lemma-inhomo-bdry-data}. Because invertibility is an open condition in the Banach space category and $E_{\mu}$ is hypo-elliptic, it is also invertible for small $\mu$. The existence and uniqueness statements follow.

To prove the tame estimate, we apply Proposition \ref{Prop-Apriori}
and get 
\begin{align}
\|\phi\|_n 
\lesssim
\|\gamma\|_n  + |\xi|_{n+2} + |\rho|_{n+1}
+
(\|\mu\|_{n+s}+1) (\|\phi\|+|x|).
\label{unisolveEst}
\end{align}
Using Lemma \ref{lemma-perpToKer} with $\mu=0$
and adding $|x|$ on both sides gives
\begin{align}
\|\phi\| +|x|
\lesssim 
\|E_0\phi\|_l+ |x| +|\xi|_{l+2} + |q_0(\phi)|_{l+1}.
\end{align} 
Because $\mathcal{H}^{p,q}$ is transversal to the image of $E_0$,
we have $\|E_0\phi\|_l +|x| \lesssim \|\gamma\|_l$.
Therefore, we get
\begin{align}
\|\phi\|+|x|
\lesssim
\|\gamma\|_l + |\xi|_{l+2} + |\rho|_{l+1}
+
\|(E_{\mu}-E_0)\phi\|_{l}
+
|(q_{\mu}-q_0)(\phi)|_{l+1}.
\label{Est1}
\end{align}
Applying 
Lemma \ref{lemma-adaptedMoserEstimate}
shows
that there exists $r$ such that
\begin{gather}
\|(E_{\mu}-E_0)\phi\|_{l}
\lesssim
\|\mu\|_r \|\phi\|_{l+2}
\nonumber
\\
\lesssim
\|\gamma\|_{l+2}+|\xi|_{l+4}+|\rho|_{l+3}
+
\|\mu\|_r |x| + \|\mu\|_r (\|\mu\|_{r+s}+1) \|\phi\|,
\label{Est2}
\end{gather}
where in the second step we have applied the a priori estimate from Proposition \ref{Prop-Apriori}.

$q_{\mu}$ is the restriction to the boundary of an operator $\hat{q}_{\mu}$ defined on all of $Z$ and analogous to above there is a t such that
\begin{gather}
|(q_{\mu}-q_0)(\phi)|_{l+1}
\lesssim
\|(\hat{q}_{\mu}-\hat{q}_0)(\phi)\|_{l+2}
\lesssim
\|\mu\|_t \|\phi\|_{l+3}
\nonumber
\\
\lesssim
\|\gamma\|_{l+3}+|\xi|_{l+5}+|\rho|_{l+4}
+
\|\mu\|_t |x| + \|\mu\|_t (\|\mu\|_{t+s}+1) \|\phi\|.
\label{Est3}
\end{gather}
If $\mu$ is sufficiently small then the estimates \eqref{Est1}, \eqref{Est2}, \eqref{Est3} after rearranging give
\begin{align}
\|\phi\|+|x|
\lesssim
\|\gamma\|_b + |\xi|_{b+2} + |\rho|_{b+1},
\label{unisolveEst2}
\end{align}
where $b=l+3$.
\eqref{unisolveEst} and
\eqref{unisolveEst2}  prove the statement.
\end{proof}

\begin{notation}
We write $G_{\mu}\gamma:= \phi$ and $H_{\mu}\gamma:=x$, where $(\phi,x)$ is the solution in Proposition \ref{Prop-UniSolveH0} with $\xi=0$ and $\rho=0$.
\end{notation}

\begin{lemma}\label{lemma-commutator}
We have
\begin{align*}
[\dbar_{\mu},G_{\mu}]\gamma = T_1(\mu)\{N_{\mu},\gamma \} + T_2(\mu)(H_{\mu}\gamma),
\end{align*}
where $T_1(\mu)$ is a map which satisfies a tame estimate in $\mu$ and is linear in the Nijenhuis tensor $N_{\mu}$ of $I_{\mu}$,
and the map $T_2(\mu)$ is $\mu$-tamely small (see Definition \ref{def-mu-small}). In particular, the commutator vanishes for $\mu=0$.
\end{lemma}
\begin{proof}
Suppose we have $G_{\mu}\gamma =\phi$, $H_{\mu}\gamma = x$, i.e.
\begin{align*}
E_{\mu}\phi + x =\gamma,
\quad
p(\phi) = 0,
\quad
q_{\mu}(\phi) = 0.
\end{align*}
Applying $\dbar_{\mu}$ we get
\begin{align*}
E_{\mu}\dbar_{\mu}\phi + [\dbar_{\mu},E_{\mu}]\phi + \dbar_{\mu}x
=
\dbar_{\mu}\gamma.
\end{align*}
Thus if we define $(\varphi_1,h_1)$ and $(\varphi_2,h_2)$ as the solutions of
\begin{align*}
E_{\mu}\varphi_1 + h_1 = [\dbar_{\mu},E_{\mu}]\phi,
\quad
p(\varphi_1) = 0,
\quad
q_{\mu}(\varphi_1) = -q_{\mu}(\dbar_{\mu}\gamma),
\end{align*}
and 
\begin{align*}
E_{\mu}\varphi_2 + h_2 = \dbar_{\mu}x,
\quad
p(\varphi_2) = 0,
\quad
q_{\mu}(\varphi_2) = 0,
\end{align*}
we get $G_{\mu}\dbar_{\mu}\gamma = \dbar_{\mu}\phi + \varphi_1 + \varphi_2$ and thus 
\begin{align*}
[G_{\mu},\dbar_{\mu}]\gamma = \varphi_1 + \varphi_2.
\end{align*}
By Proposition \ref{Prop-UniSolveH0}
$\varphi_1$ satisfies a tame estimate in $\mu$. 
Furthermore, we have
\begin{align*}
[\dbar_{\mu},E_{\mu}]
=
\dbar_{\mu}^2\dbar_{\mu}^*-\dbar_{\mu}^*\dbar_{\mu}^2
=
[\dbar_{\mu}^2,\dbar_{\mu}^*]
=
-[\partial_{\mu} i_{N''_{\mu}}+i_{N''_{\mu}}\partial_{\mu},\dbar^*_{\mu}]
\end{align*}
and
\begin{align*}
q_{\mu}(\dbar_{\mu}\phi)
=
p(\dbar_{\mu}^2\phi)
=
-p(\partial_{\mu} i_{N''_{\mu}}\phi+i_{N''_{\mu}}\partial_{\mu}\phi).
\end{align*}
Therefore, $\varphi_1$ depends linearly of $N_{\mu}$
and $T_1(\mu)\{N_{\mu},\gamma\}:= \varphi_1$ has the desired properties.

Because $\dbar_0 x =0$, we have $\dbar_{\mu}x = (\dbar_{\mu}-\dbar_{0})x$. Applying Lemma \ref{lemma-adaptedMoserEstimate} to the operator $\dbar_{\mu}-\dbar_0$, we see that $T_2(\mu)(H_{\mu}\gamma):=\varphi_2 = G_{\mu}\dbar_{\mu}H_{\mu}\gamma$ has the desired properties.
\end{proof} 

\subsection{Geometry on the boundary}

Suppose that the $6$-manifold $Z$ with boundary $M=\partial Z$
has a torsion-free $\SL(3,\bC)$-structure $\Psi+ i\hat{\Psi}$.
The restrictions of $\Psi$ and $\hat{\Psi}$ induce closed $3$-forms 
$\psi$ and $\hat{\psi}$ on $M$. In \cite{embedding-paper} the authors have shown that the boundary is strongly pseudoconvex if and only if $\psi$ and $\hat{\psi}$ satisfy certain open algebraic conditions. In this case $\psi$ and $\hat{\psi}$ will be called strongly pseudoconvex $3$-forms. We will recall their structure here.
They decompose as 
\begin{align}
\psi=\theta\wedge\alpha,
\quad
\hat{\psi}=\theta\wedge\beta,
\label{3-form-decomposition}
\end{align}
where $\theta$ is a contact $1$-form with associated contact distribution $H$ and Reeb field $v$.
This induces decompositions
\begin{align*}
TM=\bR v\oplus H, 
\quad
\Lambda^p T^*M 
=
\theta\wedge\Lambda_H^{p-1}\oplus \Lambda^p_H,
\end{align*}
where $\zeta\in\Lambda^p_H$ if $v\lrcorner\zeta=0$.
Setting $\omega:=d\theta$, the decomposition \eqref{3-form-decomposition} is normalised such that $(\omega,\alpha,\beta)$ is an orthonormal triple of sections of $\Lambda^2_H$, which means that they satisfy 
\begin{align}
\omega^2 = \alpha^2=\beta^2 
\label{squares}
\end{align}
and
\begin{align}
\omega\wedge\alpha = \omega\wedge\beta=\alpha\wedge\beta=0.
\label{wedgePairs}
\end{align}

Together with the volume form $\vol_H =\alpha^2$, this orthonormal triple defines a Euclidean structure $g_H$ on $H$, which gives the Webster metric $g=\theta^2 + g_H$ on $M$. We note that $g$ is not necessarily the restriction of a chosen background Hermitian metric $h$. 

As in $4$-dimensional Riemannian geometry we have a splitting
\begin{align}
\Lambda^2_H = \Lambda_H^+\oplus \Lambda^-_H
\label{2formSplitting}
\end{align}
of $2$-forms on $H$ into self- and anti-self-dual parts. $H$ is the real part of the complex tangent space at the boundary, and we have
\begin{align}
\Lambda^{2,0}_H = \bC (\alpha+i\beta),
\quad
\Re \Lambda^{1,1}_H = \bR \omega \oplus \Lambda_H^-.
\label{form-decomposition-boundary}
\end{align}
We will also need to consider $\SL(3,\bC)$-structures which are merely closed but not torsion-free, i.e. which satisfy $d\Psi=0$ but not $d\hat{\Psi}=0$. In this case, all of the above stays true except that now there is a non-zero $\lambda$  
such that
\begin{align}
\omega\wedge\beta = \lambda \beta^2.
\label{lambda}
\end{align}
If $\lambda<1$, then 
we have an orthonormal triple $(\tilde{\omega},\alpha,\beta)$ which spans $\Lambda^+_H$, where
\begin{align*}
\tilde{\omega}
=
\frac{1}{\sqrt{1-\lambda^2}}(\omega-\lambda\beta).
\end{align*}

In the next Lemma we relate the decomposition of forms on $Z$ induced by the complex structure $I_{\Psi}$ to the decomposition of forms on $M$ discussed above.
From the above discussion it follows that we have a direct sum decomposition $TZ|_M=\bR\, Iv \oplus \bR\, v \oplus H$. This allows us to extend sections of $\Lambda^{\bullet}T^*M$ to sections of $(\Lambda^{\bullet}T^*Z)|_M$.
\begin{lemma}\label{lemma-bdry-form-decomp}
There exists a boundary defining function $r$ such that $\theta=I dr$ and
\begin{align*}
(\Lambda^2 T^*Z)\bdry
&=
\bR\, dr\wedge\theta \oplus  dr\wedge\Lambda^1_H \oplus \theta\wedge\Lambda^1_H \oplus \Lambda^2_H,
\\
(\Re\Lambda^{1,1}_Z)\bdry
&=
\bR\, dr\wedge \theta
\oplus \{ dr\wedge\eta+\theta\wedge I\eta: \eta\in\Omega^1_H\}
\oplus \bR \tilde{\omega} \oplus \Lambda^-_H,
\\
\Lambda^2_6\bdry 
&= \{ dr\wedge\eta - \theta\wedge I\eta: \eta\in\Omega^1_H\}
\oplus \SPAN\{\alpha,\beta \}.
\end{align*}
\end{lemma}
\begin{proof}
Given any boundary defining function $r$, because $dr$ restricted to $TM$ vanishes, it must be a non-zero multiple of $I\theta$. By rescaling $r$ appropriately, we get $Idr=\theta$. The rest follows.
\end{proof}

\textbf{Notation}: If $\sigma\in\Omega^{\bullet}(Z)$, by $\sigma\bdry$ we denote the restriction of $\sigma$ to $M$ as a section of $\Lambda^{\bullet} T^*M$, and by $\sigma\|_M$ we denote the restriction of $\sigma$ to $M$ as a section of $(\Lambda^{\bullet}T^*Z)|_M$. 
By Lemma \ref{lemma-bdry-form-decomp} if $\sigma\in\Re \Omega^{1,1}$ or $\sigma\in\Omega^2_6$, then $\sigma\bdry=0$ and $\sigma\|_M=0$ are equivalent.
By $\sigma\bdryf$ we denote the component in $\mathcal{C}^{\infty}(M)\omega\oplus \Omega^-_H$ of $\sigma\|_M$. 

\begin{lemma}
\label{lemma-extension}
Let $F$ be a Hermitian vector bundle over $Z$.
\begin{compactenum}[(i)]
\item Each $s\in\Gamma(M, F|_M)$ has an extension $\mathcal{E}_F(s)\in\Gamma(Z,F)$ which satisfies the estimate
\begin{align*}
\|\mathcal{E}_F(s)\|_k \lesssim \|s\|_{k,\partial},
\end{align*}
for all $k$,
where $\| \cdot \|_{k,\partial}$ denotes the $L^2_k$-norm on the boundary.
\item Denoting the exterior derivatives on $Z$ and $M$ by $d_Z$ and $d_M$, respectively, we have
\begin{align*}
(d_Z\mathcal{E}_{\Lambda^k}\gamma)\bdry = d_M \gamma
\end{align*}
for all $\gamma\in\Gamma(M, (\Lambda^k T^*Z)|_M)$.
\end{compactenum}
\end{lemma}
\begin{proof}
(i) follows as in \cite[Theorem 1.3.7]{Schwarz-ellBVP} by choosing a collar neighbourhood of the boundary and a cut-off function in the radial direction.
\\
(ii) follows because the cut-off function is constant in a neighbourhood of the boundary.
\end{proof}

In the application of  Nash--Moser theory to prove Theorem \ref{MainThm}, we will need to consider all closed $\SL(3,\bC)$-structures $\Psi$ in a neighbourhood of a torsion-free $\SL(3,\bC)$-structure $\Psi_0$
with strongly pseudoconvex boundary.
Some of these are not torsion-free, but if the neighbourhood is chosen small enough, then $\lambda$ from \eqref{lambda} satisfies $\lambda<1$ and Lemma \ref{lemma-bdry-form-decomp} applies.

Next we make use of the analytic results from \cite{embedding-paper}.
Given a contact-metric structure $(\theta,g_H)$ defined by $\Psi_0$, 
denote by $d_H: \Omega^p_H \rightarrow \Omega^{p+1}_H$ the projection of the exterior derivative to $H$, by $d_H^*$ its formal adjoint, and $d_H^-: \Omega_H^1\rightarrow \Omega^-_H$ the anti-self-dual part. 
On $\Omega^{\bullet}_H$ the operator $\Delta_H=d_H d_H^* + d_H^* d_H$ is  an analogue of the Hodge Laplacian. In \cite{embedding-paper}, the authors have shown that $\Delta_H$ is sub-elliptic on $\Omega^-_H$ and that its kernel
\begin{align*}
\mathcal{H}_M= \{ \sigma\in\Omega^-_H: \Delta_H \sigma=0\}=
\{\sigma\in\Omega^-_H: d_H\sigma=0\}
\end{align*} 
is finite dimensional and the obstruction space to solve the equation $d_H^-\eta =\sigma$.

\begin{proposition}
\label{Prop-eliminatingBdry}
Let $\Psi_0$ be a torsion-free $\SL(3,\bC)$-structure on $Z$ with strongly pseudoconvex boundary and suppose that $\mathcal{H}_M=0$
for the contact metric structure induced by $\Psi_0$ on $M$.
Then for all closed $\Psi$ in a neighbourhood of $\Psi_0$ and all real
 $\sigma\in\Omega^{1,1}(\Psi)$ there are 
$\eta_{\Psi}(\sigma)\in\Omega^1(Z)$ and $\tau_{\Psi}(\sigma)\in\Omega^2_6(\Psi)$ such that $(\sigma-d\eta_{\Psi}(\sigma)+\tau_{\Psi}(\sigma))\bdry=0$ and the mappings $\sigma\mapsto \eta_{\Psi}(\sigma)$ and $\sigma\mapsto \tau_{\Psi}(\sigma)$ are tame.
\end{proposition}
\begin{proof}
Let $\sigma\in\Omega^2(Z)$.
By the results in \cite{embedding-paper} there exists $\eta\in\Omega^1(M)$ which depends tamely on $\Psi$ such that $\sigma\|_M - d_M\eta$ has no component in $\mathcal{C}^{\infty}(M)\omega_{\Psi} \oplus \Omega^-_H(\Psi)$, where $\omega_{\Psi}$ and $\Omega^-_H(\Psi)$ are the $\omega$ and $\Omega^-_H$ as above induced by $\Psi$. Then there exists a unique $\tau\in\Gamma(M,\Lambda^2_6(\Psi)|_M)$ such that $\sigma\|_M - d_M\eta+\tau \in \Gamma(M,dr\wedge (\Lambda^1 T^* Z)|_M)$.
The result follows from Lemma \ref{lemma-extension} if we set 
$\eta_{\Psi}(\sigma)  = \mathcal{E}_{\Lambda^1}(\eta)$
and $\tau_{\Psi}(\sigma) = \mathcal{E}_{\Lambda^2}(\tau)$.
\end{proof}

\section{Gauge fixing}

\label{section-gauge}

\subsection{Gauge fixing for the domain of the torsion-free $\SL(3,\bC)$-equation}

Denote by $\mathcal{G}$ the component of the identity of the group of diffeomorphisms
of $Z$ which fix the boundary point-wise. $\mathcal{G}$ acts on the space of stable $3$-forms, and on a given enhancement class. We want to study the torsion-free $\SL(3,\bC)$-equation modulo the action of $\mathcal{G}$. 
The tangent space of $\mathcal{G}$ at the identity is the space $\mathcal{V}_0$ of vector fields which vanish on the boundary. 
The tangent space in $\Omega^3_s(Z)$ to the $\mathcal{G}$-orbit of a given 
stable $3$-form $\Psi_0$  consists of all Lie derivatives $\Lie_X\Psi_0$
with $X\in \mathcal{V}_0$.
If $\Psi_0$ is closed, then
\begin{align*}
\Lie_X \Psi_0 = d(i_X\Psi_0).
\end{align*}
By \eqref{2-forms}
$2$-forms of the form $\tau := i_X\Psi_0$ are precisely the sections of $\Lambda^2_6$. Because $X$ vanishes on the boundary, we have $\tau\|_M=0$, or equivalently, $\tau\bdry=0$. 
Therefore, we have 
\begin{align*}
T_{\Psi_0}(\mathcal{G}\cdot \Psi_0) 
=
d\{\tau\in\Omega^2_6: \tau\bdry =0 \}.
\end{align*}
The standard approach to construct a slice for the diffeomorphism action is  to consider $3$-forms which are $L^2$-orthogonal to $T_{\Psi_0}(\mathcal{G}\cdot \Psi_0)$.
If $\tau\in\Omega^2_6$ satisfies $\tau\bdry=0$ and $\rho\in\Omega^3(Z)$, we can integrate by parts to get
\begin{align*}
(d\tau,\rho)=(\tau,d^*\rho) = (\tau,\pi_6 d^*\rho).
\end{align*}
Thus the $L^2$-orthogonal complement of $T_{\Psi_0}(\mathcal{G}\cdot \Psi_0) $ in $\Omega^3(Z)$ is the slice
\begin{align*}
S = \{ \rho \in \Omega^3(Z): \pi_6 d^*\rho =0 \}.
\end{align*}
We will prove the following slice theorem.
\begin{proposition}
\label{prop-gauge}
There is a neighborhood $U$ of $\Psi_0$ in $\Omega_s^3(Z)$ such that if $\Psi\in U$, there are neighborhoods $U_{\mathcal{G}}$ of $\mathbbm{1}$ in $\mathcal{G}$, $U_S$ of $0$ in $S$ and $U_{\Omega}$ of $\Psi$ in $\Omega^3(Z)$, such that for each $\Psi'\in U_{\Omega}$ there are unique $\varphi \in U_{\mathcal{G}}$ and $\mu \in U_S$ such that $\Psi' =\varphi^* (\Psi + \mu)$.
\end{proposition}
To prove Proposition \ref{prop-gauge}, let $\Psi$ be close to $\Psi_0$. We want to apply the Nash--Moser inverse function theorem
to the map
\begin{align*}
\mathcal{M}_{\Psi}:\mathcal{G}\times S \rightarrow \Omega^3(Z),
\quad
(\varphi, \mu)
\mapsto \varphi^*(\Psi+\mu).
\end{align*}
Given $\varphi\in\mathcal{G}$, by sending $X\in\mathcal{V}_0$ to the path $\exp(s X)\circ \varphi$, we can identify $T_{\varphi}\mathcal{G}$ with  $\mathcal{V}_0$.
The derivative of $\mathcal{M}_{\Psi}$ at $(\varphi,\mu)$ is the map
\begin{align*}
D\mathcal{M}_{\Psi}(\varphi,\mu):
\mathcal{V}_0\oplus S \rightarrow \Omega^3(Z),
\quad
(X,\rho)
\mapsto
\varphi^*(\mathcal{L}_X\Psi + \Lie_X\mu + \rho).
\end{align*}
Thus to invert the map $D\mathcal{M}_{\Psi}(\varphi,\mu)$ we need to solve the equation
\begin{align}
\mathcal{L}_X \Psi_0 + \rho
+  
 \mathcal{L}_X(\psi+\mu)
=
(\varphi^{-1})^*\gamma ,
\label{linGauge}
\end{align}
where $\Psi=\Psi_0+\psi$.
We will first treat the case $\psi=0$ and $\mu=0$.
To decompose a given $3$-form $\gamma$
as $\Lie_X\Psi_0+\rho=\gamma$, write $\rho=\gamma+\chi$.
This means that we need to find a $3$-form $\chi$ such that 
$\pi_6 d^*(\gamma+\chi)=0$, i.e. we need to solve the equation $\pi_6 d^*\chi = -\pi_6 d^*\gamma$. We can use Lemma \ref{lemma-6BVP}
and take $\chi = -dG_6(\pi_6 d^*\gamma)$.
Denote by 
\begin{align*}
V: \Lambda^2_6\rightarrow TZ
\end{align*}
the bundle isomorphism which is the inverse of $X\mapsto i_X\Psi_0$.
Thus we get a map
\begin{align*}
W=(W_1,W_2): \Omega^3(Z) \rightarrow \mathcal{V}_0 \oplus S,
\quad
\gamma
\mapsto (V(G_6(\pi_6 d^*\gamma)),\gamma - d G_6(\pi_6 d^*\gamma)),
\end{align*}
which is characterised by
\begin{align*}
\mathcal{L}_{W_1(\gamma)}\Psi_0 + W_2(\gamma) = \gamma.
\end{align*}
Because of the ellipticity of the boundary value problem in Lemma \ref{lemma-6BVP}, there are positive constants $C_n'$ such that
\begin{align}
\|W_1(\gamma)\|_n \leq C_n' \|\gamma\|_{n-1},
\quad
\|W_2(\gamma)\|_n \leq C_n' \|\gamma\|_n.
\label{constants}
\end{align}

\begin{lemma}
The map $W$ is an isomorphism.
\end{lemma}
\begin{proof}
We first prove that $W$ is injective. Suppose that $W(\gamma)=0$. Then $G_6(\pi_6 d^*\gamma)=0$ and thus $\pi_6 d^*\gamma=0$. Because $W_2(\gamma)=0$, we get $\gamma=0$. 

$W$ clearly is the identity on $S\subset \Omega^3(Z)$, and given $X\in\mathcal{V}_0$, by Lemma \ref{lemma-6BVP} $W$ maps $d(i_X\Psi)$ to $(X,0)$. Thus $W$ is also surjective.
\end{proof}

Now we are ready to solve equation \eqref{linGauge} for small non-zero $\psi$.
\begin{lemma}
\label{lemma-gaugeEqu}
There exists $\epsilon > 0$, such that if $\kappa\in\Omega^3(Z)$ satisfies $[[\kappa]]_1 < \epsilon$, then given $\gamma\in\Omega^3(Z)$ there is a solution $W(\kappa)\gamma=(W_1(\kappa)\gamma,W_2(\kappa)\gamma):=(X,\rho)\in \mathcal{V}_0\oplus S$ of the equation
\begin{align}
\mathcal{L}_X \Psi_0 + \rho 
+ 
\mathcal{L}_X\kappa
=
\gamma,
\label{LieEquation}
\end{align}
and the solution satisfies the tame estimate
\begin{align}
\|W(\kappa)\gamma\|_n
\lesssim
\|\gamma\|_n + [[\kappa]]_{n+1} \|\gamma\|.
\label{Wtame}
\end{align}
\end{lemma}
\begin{proof}
On $\Omega^3(Z)$ define the operator $L(\kappa)\xi = - \Lie_{W_1(\xi)}\kappa$. To solve \eqref{LieEquation}, we need to solve the equation 
$(\mathbbm{1}-L(\kappa))\xi = \gamma$. 
Indeed, if $W(\xi)=(X,\rho)$, i.e. $\xi = \Lie_X\Psi_0 +\rho$, then
$(\mathbbm{1}-L(\kappa))\xi = \Lie_X\Psi_0 +\rho + \Lie_X\kappa$.

By Lemma \ref{lemma-adaptedMoserEstimate} we have
\begin{align*}
\|L(\kappa)\xi\|_n
\leq (sc) [[\kappa]]_1 \|W_1(\xi)\|_{n+1}
+
(lc)
[[\kappa]]_{n+1}
\|W_1(\xi)\|_1,
\end{align*}
where (sc) denotes an arbitrarily small constant, and (lc) a possibly large one. If we choose (sc) less than $1/C'_{n+1}$, we thus get constants $C_n$ such that
\begin{align*}
\|L(\kappa)\xi\|_n
\leq  [[\kappa]]_1 \|\xi\|_n
+
C_n
[[\kappa]]_{n+1}
\|\xi\|.
\end{align*}
By Proposition \ref{prop-NeumannSeries}, $\mathbbm{1}-L(\kappa)$ has a tame inverse if $\kappa$ is sufficiently small.
\end{proof}
Suppose that $[[\psi]]_1 < \epsilon/2$ and $[[\mu]]_1 < \epsilon/2$.
Set 
\begin{align*}
W^{\Psi}(\mu)=(W^{\Psi}_1(\mu),W^{\Psi}_2(\mu)):= W(\psi+\mu).
\end{align*}
$(X,\rho) = (W_1^{\Psi}(\mu)\gamma,W_2^{\Psi}(\mu)\gamma)$ is characterised by
\begin{align}
\label{def-W}
\mathcal{L}_X\Psi + \rho + \mathcal{L}_X\mu = \gamma.
\end{align}
By Lemma \ref{lemma-gaugeEqu} $W^{\Psi}(\mu)\circ (\varphi^{-1})^*$
is an inverse of $D\mathcal{M}_{\Psi}(\varphi,\mu)$.
By \eqref{Wtame} we have
\begin{align*}
\|W^{\Psi}(\mu)\gamma\|_n
\lesssim \|\gamma\|_n + (1 + [[\mu]]_{n+1}) \|\gamma\|.
\end{align*}
Therefore, $W^{\Psi}(\mu)$ is tame in $\mu$, and $W^{\Psi}(\mu)\circ (\varphi^{-1})^*$ is tame as a composition of smooth tame maps.
This completes the proof of Proposition \ref{prop-gauge}.

Next we want to consider a slice for the action of $\mathcal{G}$ in an enhancement class. For this we define
\begin{align*}
T=S\cap \mathcal{E}^2_D
=
\{ d\sigma: \sigma\in\Omega^2(Z), \sigma\bdry=0, \pi_6 d^*d\sigma=0\},
\end{align*}
where $\mathcal{E}^2_D$ is the space from \eqref{exact-forms}.
Then Proposition \ref{prop-gauge} implies

\begin{corollary}
\label{lemma-nearbyClass}
Suppose that $\Psi \in U$ and $d\Psi=0$. Then there is a neighborhood $U_{\Psi}$ of $\Psi$ in its enhancement class, a neighbourhood $U_{\mathcal{G}}$ of $\mathbbm{1}$ in $\mathcal{G}$ and a neighbourhood $U_T$ of $0$ in $T$, such that for each $\Psi' \in U_{\Psi}$ there are unique $\varphi\in U_{\mathcal{G}}$ and $\rho\in U_T$ such that $\Psi' = 
\varphi^*(\Psi + \rho)$.
\end{corollary}
We conclude this section by identifying the space $T$ with certain $2$-forms.
Write
\begin{align*}
A'=\{\sigma\in\Omega^{1,1}: d^*\sigma=0,\sigma\|_{\partial,4}=0 \}.
\end{align*}
Next we construct a linear map
\begin{align*}
j: \mathcal{H}^2_D\rightarrow A'.
\end{align*}
If $h\in\mathcal{H}^2_D$, then there is $1$-form $\eta$ which solves   the boundary value problem \eqref{ellBVP11} with $\lambda = -d^* h_{1,1}$. Then $h_{1,1}+d_{1,1}\eta$ is coclosed. Because $h\bdry =0$ and $\eta\bdry=0$, we have $(h+d\eta)\bdry =0$. Therefore by Lemma \ref{lemma-bdry-form-decomp} we have $(h_{1,1}+d_{1,1}\eta)\bdryf =0$,
and thus 
$h_{1,1} + d_{1,1}\eta \in A'$. 
By Lemma \ref{lemma-EllBdry}
\begin{align}
j(h):= h_{1,1} + d_{1,1}\eta \in A'
\label{def-j}
\end{align} 
is well-defined.

Define $A$ to be the $L^2$-orthogonal complement of $\im( j)$ in $A'$.
We now define a linear  map $A\rightarrow T$ as follows.
By Lemma \ref{lemma-bdry-form-decomp} there is a bundle map $\tau: \Omega^2_{1,1}(Z)\rightarrow \Omega^2_6(Z)$ supported in a neighbourhood of the boundary  such that if $\zeta\in\Omega^{1,1}(Z)$ satisfies $\zeta\bdryf=0$, we have $(\zeta+\tau(\zeta))\bdry=0$. 
If $\sigma\in A'$, set 
\begin{align*}
\Gamma(\sigma)
=
\tau(\sigma)-G_6(\pi_6 d^*d(\sigma+\tau(\sigma)).
\end{align*}
Then by construction $d(\sigma+\Gamma(\sigma))\in T$.
\begin{lemma}
\label{lemma-2forms}
The map $d\circ (\mathbbm{1}+\Gamma):A'\rightarrow T$ is surjective and has kernel $\im(j)$. Thus $d\circ(\mathbbm{1}+\Gamma)$ restricted to $A$ is an  isomorphism. 
\end{lemma}
\begin{proof}
We first show that $d\circ(\mathbbm{1}+\Gamma)$ is surjective.
Suppose that $d\zeta\in T$. By Lemma \ref{lemma-EllBdry} there exists $\eta\in\Omega^1(Z)$ with $d^*d_{1,1}\eta=-d^*\zeta_{1,1}$ and $\eta\bdry=0$. Therefore, by replacing $\zeta$ with $\zeta-d\eta$, we can assume that $\sigma:=\zeta_{1,1}\in A'$. Then we have
\begin{align*}
\zeta = \sigma + \tau(\sigma) + \xi, 
\end{align*}
where $\xi:= \zeta_6-\tau(\sigma)$. Because $\zeta\bdry=0$ and $(\sigma+\tau(\sigma))\bdry=0$, we have $\xi\bdry=0$, and thus by Lemma \ref{lemma-bdry-form-decomp} $\xi\bdrytn =0$. 
By the definition of $T$, we have  $\pi_6 d^*d\zeta=0$ and thus
$\pi_6 d^*d\xi = - \pi_6 d^*d(\sigma+\tau(\sigma))$.
By the uniqueness statement in Lemma \ref{lemma-6BVP}  we must have $\xi = - G_6(\pi_6 d^*d(\sigma+\tau(\sigma))$,
i.e. $d\zeta = d(\mathbbm{1}+\Gamma)(\sigma)$.

Next we show $\ker (d\circ(\mathbbm{1}+\Gamma)) \subset \im(j)$.
Suppose we have  $d(\sigma+\Gamma(\sigma))=0$. By Proposition \ref{Hodge-mixed} (i)
there are $h\in\mathcal{H}^2$ and $\eta\in\Omega^1(Z)$ with $\eta\bdry=0$ such that 
\begin{align*}
\sigma+\Gamma(\sigma)
= 
h + d\eta.
\end{align*} 
Because $(\sigma+\Gamma(\sigma))\bdry=0$ and $\eta\bdry=0$, we have $h\in\mathcal{H}^2_D$.
Projecting onto the component of type $(1,1)$ gives $\sigma= h_{1,1} + d_{1,1}\eta$. Because $d^*\sigma=0$, $\eta$ is a solution of the boundary value problem 
from Lemma \ref{lemma-EllBdry} 
with $\lambda = - d^*h_{1,1}$. Therefore $\sigma=j(h)$.

Next we show $\im(j) \subset \ker (d\circ(\mathbbm{1}+\Gamma))$.
Let $\sigma = h_{1,1}+d_{1,1}\eta \in\im(j)$ as defined in \eqref{def-j}. If we write
\begin{align*}
h + d\eta = \sigma+\tau(\sigma) + \xi,
\end{align*}
then as above we can conclude $\xi\bdrytn =0$. Because $\pi_6 d^*d$ of the left-hand side vanishes, we must have $\xi = - G_6(\pi_6 d^*d(\sigma+\tau(\sigma))$, and thus $h + d\eta = \sigma + \Gamma(\sigma)$. Thus $d(\sigma+\Gamma\sigma)=0$.
\end{proof}

Finally we characterise $A$ to be the set of co-exact forms in $A'$.
\begin{lemma}
\label{lemma-charA}
We have
\begin{align*}
A=\{\sigma\in\Omega^{1,1}: \sigma\, \text{is co-exact}, \sigma\bdryf=0\}.
\end{align*}
\end{lemma}
\begin{proof}
We first show that $\sigma\in A$ is co-exact.
By Proposition \ref{HodgeDirichlet} (i) we can write $\sigma = h + d^*\mu$ with $h\in\mathcal{H}^2_D$ and $\mu\bdry =0$ because $\sigma$ is co-closed.
By the definition of $A$, $\sigma$ is $L^2$-orthogonal to $j(h)\in A'$. Thus we have
\begin{align*}
0 = (\sigma,j(h)) = (h+d^*\mu,h_{1,1}+d_{1,1}\eta) = (h+d^*\mu,h+d\eta) = \|h\|^2.
\end{align*}
This implies $h=0$ and that $\sigma$ is co-exact.

Vice versa, let $\sigma \in A'$ be co-exact with $\sigma=d^*\mu$.
Then for all $h\in\mathcal{H}^2_D$ we have
\begin{align*}
(\sigma,j(h)) = (d^*\mu, h_{1,1}+d_{1,1}\eta)
= (d^*\mu,h+d\eta) = 0.
\end{align*}
\end{proof}

\subsection{Gauge fixing for target space of the torsion-free $\SL(3,\bC)$-equation}

Denote by
\begin{align*}
\mathcal{X}=\{ \Psi\in\Omega^3_s(Z): d\Psi=0\}
\end{align*}
the space of closed $\SL(3,\bC)$-structures.
$\Psi\in\mathcal{X}$ defines a torsion-free $\SL(3,\bC)$-structure if and only if it is a zero of the map $\Psi\in \mathcal{X}\mapsto dP(\Psi)$. It is convenient to study instead a map which takes values in $2$-forms. Using the Hodge star operator of a Hermitian metric for our reference torsion-free $\SL(3,\bC)$-structure $\Psi_0$,
we set $F(\Psi):= *dP(\Psi)$. By Lemma \ref{lemma-22} the torsion $dP(\Psi)$ has type $(2,2)$ with respect to the almost-complex structure defined by $\Psi$. Thus we can regard $F$  a section $F:\mathcal{X}\rightarrow \mathscr{C}$ of the  vector bundle
\begin{align*}
\mathscr{C}=\coprod_{\Psi\in\mathcal{X}} C_{\Psi}
\rightarrow \mathcal{X},
\end{align*}
where the fiber over $\Psi$ is
\begin{align*}
C_{\Psi}=\{\zeta\in\Omega^2(Z): \zeta\, \text{is co-exact}, \zeta\in\Gamma(*\Lambda^{2,2}_{\Psi}) \}.
\end{align*}
We want to apply the Nash--Moser inverse function theorem and need to make sense of a derivative of $F$. To this end 
we will construct a tame trivilisation $R: \mathscr{C}|_U \rightarrow U\times C$ for $\mathscr{C}$ over a neighbourhood $U$ of $\Psi_0$, where we write
\begin{align*}
C=C_{\Psi_0}=\{\zeta\in\Omega^{1,1}: \zeta\, \text{is co-exact}\}
.
\end{align*}

\begin{remark}
\label{remark-A<C}
By Lemma \ref{lemma-charA}  we have $A\subset C$.
\end{remark}

Given a co-exact $2$-form $\sigma = d^*\xi$, decompose $\sigma = \sigma_{1,1}+\sigma_6$, where $\sigma$ is a real form of type $(1,1)$ and $\sigma_6\in\Omega^2_6(Z)$, and solve the equation $\pi_6d^*d\zeta = \sigma_6$ with Lemma \ref{lemma-6BVP}. Then $
\sigma-d^*dG_6\sigma_6 = d^*(\xi-dG_6\sigma_6)
\in C$. 
The restriction of the mapping 
\begin{align}
\label{def-phi}
R:\sigma\mapsto \sigma-d^*dG_6\sigma_6
\end{align}
to $C_{\Psi}$
defines a family of maps
\begin{align*}
R|_{C_{\Psi}}: C_{\Psi}\rightarrow C.
\end{align*}

\begin{lemma}
There exists a neighborhood $U$ of $\Psi_0$ in $\Omega^3_s(Z)$ such that for all $\Psi\in U$ the map $R|_{C_{\Psi}}$ is bijective
and its inverse $R^{-1}(\Psi)$ satisfies for each $n$ the estimate
\begin{align}
\|R^{-1}(\Psi)\zeta\|_{n}
\lesssim \|\zeta\|_{n}+ [[\Psi-\Psi_0]]_n \|\zeta\|
\label{inverseTame}
\end{align}
for all $\zeta\in C$ and $\Psi\in U$, i.e. the map $R^{-1}$ is tame.
\end{lemma}
\begin{proof}
For $\Psi$ close to $\Psi_0$, the bundle $*\Lambda^{2,2}_{\Psi}\subset \Lambda^2 T^*Z$ is the graph of a bundle map 
\begin{align*}
\mu(\Psi):\Lambda^{1,1}\rightarrow \Lambda^2_6
\end{align*}
If we define a map
\begin{align*}
\Pi(\Psi)\zeta = \mu(\Psi)\pi_{1,1}\zeta - \pi_6\zeta,
\end{align*}
then we have the characterisation $*\Lambda^{2,2}_{\Psi} = \ker \Pi(\Psi)$. 

On $\Omega^2(Z)$ define $L(\Psi)\sigma=d^*d G_6 \mu(\Psi)\pi_{1,1}\sigma$, so that 
$R|_{C_{\Psi}}$ is the restriction of the tame map $\mathbbm{1}-L(\Psi):\Omega^2(Z)\rightarrow \Omega^2(Z)$ to $C_{\Psi}$.

By \eqref{constant6} there are constants $C_n'$ such that
\begin{align*}
\|L(\Psi)\zeta\|_n \leq C'_n\, \|\mu(\Psi)\pi_{1,1}\zeta\|_n,
\end{align*}
Thus if we choose $\epsilon$ in Lemma \ref{lemma-adaptedMoserEstimate} such that $\epsilon C'_n<1$, we get constants $C_n >0$ such that for each $n$ we have
\begin{align*}
\|L(\Psi)\zeta\|_n \leq [[\Psi-\Psi_0]]_0 \, \|\zeta\|_n + C_n [[\Psi-\Psi_0]]_n\, \|\zeta\|_0.
\end{align*}
By Proposition \ref{prop-NeumannSeries} the map $\mathbbm{1}-L(\Psi)$ is a tame isomorphism for $\Psi$ in a neighbourhood of $\Psi_0$. 
We are left to prove that the inverse maps $C$ to $C_{\Psi}$.
Let $\zeta\in\Omega^{1,1}(Z)$.
We have
\begin{align*}
\Pi(\Psi)L(\Psi)\gamma
=
\mu(\Psi)\pi_{1,1} L(\Psi)\gamma - \pi_6 L(\Psi)\gamma
=
\mu(\Psi)\pi_{1,1}(L(\Psi)\gamma-\gamma),
\end{align*}
and thus 
\begin{align*}
\Pi(\Psi)\sum_{j=0}^k L^j(\Psi)\zeta
=
\mu(\Psi)\pi_{1,1}L^k(\Psi)\zeta.
\end{align*}
Because $L^k(\Psi)\zeta\rightarrow 0$, we see that $\Pi(\Psi)\sum_{k=0}^{\infty} L^k(\Psi)\zeta=0$. Thus the limit of the Neumann series is a section of $*\Lambda^{2,2}_{\Psi}$. 
If $\zeta$ is co-exact, then  all finite sums $\sum_{j=0}^k L^j(\Psi)\zeta$ are $L^2$-orthogonal to $\mathcal{H}^2\oplus \mathcal{E}^2$. Because the series converges in $L^2$, we see from Proposition \eqref{HodgeDirichlet} (ii) that the limit must be co-exact as well. Thus $(\mathbbm{1}-L(\Psi))^{-1}$ maps $C_0$ to $C_{\Psi}$.
\end{proof}

\section{Proof of Theorem \ref{MainThm}}

\label{section-proof}

Denote by $B$ a neighbourhood of $0$ in the space of closed $3$-forms on $Z$. If $B$ is small enough, then for any $b\in B$, the $3$-form $\Psi_0+b$ is stable.
Thus $B\rightarrow\mathcal{X}, b\mapsto \Psi=\Psi_0+b$ gives a parametrisation of a neighborhood of $\Psi_0$ in the space $\mathcal{X}$ of closed $\SL(3,\bC)$-structures.
By choosing $B$ small enough, we can apply
Corollary \ref{lemma-nearbyClass} and Lemma \ref{lemma-2forms}
so that up to the action of $\mathcal{G}$ all elements in a neighbourhood of $\Psi$ in its enhancement class are of the form $\Psi+d(a+\Gamma a)$ with $a\in A$.
Therefore, if we define the map
\begin{align*}
\mathcal{F}: U\subset(A\times B) \rightarrow C,
\quad
(a,b)
\mapsto
R(*dP\Psi'),
\,
\textrm{where}\,
\Psi'=\Psi_0 + b + d(a+\Gamma a),
\end{align*}
where $U$ is a small neighbourhood of $(0,0)$ and $R$ is the map from \eqref{def-phi},
then
torsion-free $\SL(3,\bC)$-structures in the enhancement class of $\Psi=\Psi_0+b$ correspond to $a\in A$ such that $\mathcal{F}(a,b)=0$.
Proving  Theorem \ref{MainThm} amounts to showing that for each  sufficiently small $b\in B$ there is a unique $a\in A$ close to $0$ such that $\mathcal{F}(a,b)=0$.

\subsection{Existence under stronger hypothesis}

For the existence part  we will apply a variation of the Nash--Moser implicit function theorem which is due to Zehnder \cite{ZehnderIFT}. To use the classical Nash--Moser implicit function theorem,
we would have to find a tame right inverse  
\begin{align}
V: U'\times C\rightarrow  A,
\label{map-V}
\end{align}
for the partial derivative $D_{a}\mathcal{F}: U'\times A\rightarrow C$
in some neighbourhood $U'\subset U$ of zero.
We use the notation
\begin{align*}
L_{a,b}:= D_{a}\mathcal{F}(a,b):A\rightarrow C.
\end{align*}

In our situation, we will only be able to find a genuine right inverse for $L_{a,b}$
if the $\SL(3,\bC)$-structure $\Psi'$ is torsion-free, i.e. if $\mathcal{F}(a,b)=0$. If $\Psi'$ is not torsion-free, we need to allow a small error. 

The following Theorem is the core result of this paper. It leads to a proof of the existence part of Theorem \ref{MainThm} under slightly stronger hypotheses. In the next subsection we give the full proof of Theorem \ref{MainThm}.

\begin{theorem}
\label{Thm-specialCase}
Suppose that $\mathcal{H}_M=0$ and $\mathcal{H}^{1,2}=0$. Then 
there exists an open neighbourhood $U'\subset U$ of $(0,0)$, a tame map $V$ as in \eqref{map-V}, a tame map
\begin{align*}
Q: U' \times \Omega^2(Z,TZ) \times C \rightarrow C,
\end{align*}
which is linear in the $\Omega^2(Z,TZ)$-variable and linear in the $C$ variable, 
and 
a tame map
\begin{align*}
\tilde{Q}: U' \times \Omega^2(Z,TZ) \times C \rightarrow C,
\end{align*}
which is quadratic in the $\Omega^2(Z,TZ)$-variable and linear in the $C$-variable, such that for all $(a,b)\in U'$ and for all $c\in C$ we have
\begin{align*}
L_{a,b}V_{a,b}c
=
c + Q_{a,b}(N_{a,b},c)
+
\tilde{Q}_{a,b}(N_{a,b},c)
,
\end{align*}
where $N_{a,b}$ denotes the Nijenhuis tensor associated with $\Psi'$.
\end{theorem}

Given Theorem \ref{Thm-specialCase}, we now explain how the existence statement in Theorem \ref{MainThm} follows from Theorem \ref{Nash-Moser-Zehnder-Thm}. The uniqueness statement in Theorem \ref{MainThm} and the more general case $\mathcal{K}=0$ will be covered in the next section. 

By possibly shrinking $U$, we can ssume without loss of generality that we have a product $U=U_A\times U_B$, where $U_A, U_B$ are nighbourhoods of $0$ in $A,B$. $U_B$ can be chosen small enough such that the first set of hypothesis in Theorem \ref{Nash-Moser-Zehnder-Thm} are satisfied for $\phi = \mathcal{F}(\,\cdot\, ,b)$ and $u_0=0$ for all $b\in U_B$. Because $V$ is tame $\psi(a) = V_{a,b}$ satisfies the hypothesis of Theorem \ref{Nash-Moser-Zehnder-Thm}.
By Lemma \ref{lemma-N}  $N_{a,b}$ is a function of $\mathcal{F}(a,b)$ and $(a,b)$. Therefore, we can define $\mathcal{E}(a)$ by the equation
\begin{align*}
\mathcal{E}(a)(\mathcal{F}(a,b),c) = Q_{a,b}(N_{a,b},c)+\tilde{Q}_{a,b}(N_{a,b},c).
\end{align*}
By the third Moser estimate \cite[p. 440]{Hamilton1977Deformation2} we have $\|N_{a,b}\|_k \lesssim \|\mathcal{F}(a,b)\|_{k+d}$ for all $k$ and some fixed $d$.
By the structure of $Q$ and $\tilde{Q}$ there are $r,s$ such that for each $k$ we have
\begin{align*}
\|\mathcal{E}(a)(f,c)\|_k
\lesssim
\|f\|_{k+r} \|c\|_{k+s}.
\end{align*} 
Thus
all hypothesis of Theorem \ref{Nash-Moser-Zehnder-Thm} are satisfied.
Because $\mathcal{F}$ is continuous, we have $\|\phi(u_0)\|_{2d} = \|\mathcal{F}(0,b)\|_{2d}\rightarrow 0$ as $b\rightarrow 0$. Thus 
the existence 
of $a\in U_a$ with $\mathcal{F}(a,b)=0$ follows for all $b$ in some neighbourhood of $0$.

In the remainder of this subsection we will prove Theorem \ref{Thm-specialCase}.
The non-linearity of $\mathcal{F}$ is purely contained in the map $P$. 
By \cite[Propositions 4 and 5]{Hitchin}
the derivative of $P$ at $\Psi'$ is the bundle map
\begin{align*}
J_{a,b}:\Omega^3(Z)\rightarrow \Omega^3(Z)
\end{align*}
which acts acts by multiplication with $-i$ on $3$-forms which are of type $(3,0)$ and $(2,1)$ with respect to $\Psi'$, and by multiplication with $i$ on $3$-forms of type $(1,2)$ and $(0,3)$ with respect to $\Psi'$.
Thus we get
\begin{align*}
L_{a,b}= R \circ * \circ  d\circ J_{a,b}\circ d \circ (\mathbbm{1}+\Gamma).
\end{align*}

In the remainder of the proof we will refer to error terms which
have the properties of $Q$ or $\tilde{Q}$, i.e. are linear in $c$ and either linear or quadratic in $N_{a,b}$,  as \textit{small error} terms.  
\begin{lemma}
\label{lemma-6smallerror}
$dJ_{a,b}d$ applied to $\Omega^2_6(\Psi')$ only contributes small error terms.
\end{lemma}
\begin{proof}
Let $\varphi\in\Omega^{2,0}(\Psi')$.
Dropping the subscript $(a,b)$ and using the identities from Lemma \ref{lemma-ACIDs}, we compute
\begin{align*}
\pi_{3,1}dJd\varphi
&=
-\sqrt{-1}
(\partial\dbar\varphi+\dbar\partial\varphi-i_{N'}i_{N''}\varphi)
=
2\sqrt{-1} i_{N'}i_{N''}\varphi,
\\
\pi_{2,2}dJd\varphi
&=
-\sqrt{-1}
(\dbar^2\varphi+i_{N''}\partial\varphi-\partial i_{N''}\varphi)
=
2\sqrt{-1} \partial  i_{N''}\varphi,
\\
\pi_{1,3}dJd\varphi
&=
\sqrt{-1} 
(\dbar i_{N''}\varphi - i_{N''}\dbar\varphi).
\end{align*}
\end{proof}

Write
\begin{align*}
\Omega^2_D = \{\kappa \in\Omega^2(Z): \kappa\bdry=0 \}
\end{align*}
for $2$-forms which satisfy the Dirichlet boundary condition. 
To prove Theorem \ref{Thm-specialCase}, we first simplify our problem and observe that it is enough to find a tame right $\tilde{V}$ inverse up to small error terms for the family of maps
\begin{align*}
\tilde{L}_{a,b}: \Omega^2_D \rightarrow C,
\quad
\chi \mapsto R(*dJ_{a,b}d\chi).
\end{align*}
By \eqref{def-W} $\chi\in\Omega^2_D$ determines $X\in\mathcal{V}_0$
and $\rho\in T$ such that
\begin{align*}
d\chi = d(i_X\Psi')+\rho.
\end{align*}
By Lemma \ref{lemma-2forms} there exists $\alpha\in A$ such that $\rho = d(\mathbbm{1}+\Gamma)\alpha$. Because $dJ_{a,b}d$ applied to $i_X\Psi'$ is a small error term by Lemma \ref{lemma-6smallerror} we have
\begin{align*}
\tilde{L}_{a,b}\chi = L_{a,b}\alpha + \text{small error},
\end{align*}
and $\alpha$ depends tamely on $\chi$ as the map $W^{\Psi'}_2(0)$ from \eqref{def-W} is tame.

As in Section \ref{section-dbar} write $\mu = \mu(a,b)$ for the form in $\Omega^{0,1}(Z,T^{1,0}Z)$ 
which determines the almost complex structure $I_{\Psi'}$ relative to $I_{\Psi_0}$. Given $c\in C$, we will find a genuine solution of  $\tilde{L}_{a,b}\chi = c$ only
if $\Psi'$ is torsion-free, which we now explain.
There are several simplifications in this case.
At such $\Psi'$ 
the exterior derivative splits as $d=\partial_{[\mu]}+\dbar_{[\mu]}$, we have $\dbar_{[\mu]}^2=0$ and $dJ_{\Psi'}d$ vanishes on $\Omega^2_6(\mu)$ because of the diffeomorphism invariance.
  Thus we have $dJ_{a,b}d = 2i \partial_{[\mu]}\dbar_{[\mu]}\pi^{1,1}_{\mu}$ and we have to prove a variation of a 
$\partial\dbar$-Lemma. We do not assume that our Hermitian metric $h$ is K\"ahler, because it
is not naturally associated with any of our $\SL(3,\bC)$-structures. 
Instead we reduce our problem to a classical $\dbar$-equation by using  Proposition \ref{prop-top-degree-Neumann}
which naturally comes from the existence of the torsion-free $\SL(3,\bC)$-structure $\Psi_0$.

First consider $(a,b)=(0,0)$. Given $c\in C$, set $\rho:= *c$,
which is an exact $4$-form of type $(2,2)$. Write $\rho=d\beta$.
We can make the choice of the primitive $\beta$ unique by projecting on $\mathcal{C}^3_N$ in the Hodge decomposition from Proposition \ref{Hodge-mixed}.
Under the type decomposition we have $\beta = \beta_{3,0} + \beta_{2,1} + \beta_{1,2} + \beta_{0 ,3}$ with $\beta_{p,q}\in\Omega^{p,q}(Z)$ and $\bar{\beta}_{p,q}=\beta_{q,p}$. We want to eliminate the components of type $(3,0)$ and $(0,3)$. By Proposition \ref{prop-top-degree-Neumann} we can write $\beta_{0,3}=\dbar \lambda$.
Setting $\gamma=\pi_{1,2}(\beta-d(\lambda+\bar{\lambda}))$, we have $\rho = \partial \gamma + \dbar \bar{\gamma}$ and $\dbar\gamma=0$.
Using Proposition \ref{prop-dbar-Neumann}, under the condition $\mathcal{H}^{1,2}=0$ we can write $\gamma = \dbar \zeta$.
Setting $\chi=\Im\zeta$, we have $\rho=2i\partial\dbar\chi=dJd\chi$. 
We do not necessarily have  $\chi\bdry=0$,
but by Proposition \ref{Prop-eliminatingBdry}
we can adapt $\chi$ to $\tilde{\chi}\in\Omega^2_D$
by adding a $2$-form which is exact and  a $2$-form of type $6$, both of which are annihilated by $dJd$. The same argumentation goes through for $(a,b)$ in some neighborhood  of $(0,0)$
if $\Psi'$ is torsion-free. If $\mathcal{F}(a,b)\neq 0$ this will not work, but in a neighborhood of zero we can still define
\begin{align}
\rho &= *R^{-1}(\Psi')c = d\beta, \,\,(\beta\in\mathcal{C}^3_N)
\nonumber
\\
\lambda &= c_{\mu} \dbar_{\mu}^*G_{\mu}c_{\mu}^{-1} \pi^{0,3}_{\mu}\beta,
\nonumber
\\
\gamma
&=
\pi^{1,2}_{\mu}(\beta-d(\lambda+\bar{\lambda})),
\label{2step-process}
\\
\zeta &= c_{\mu}\dbar_{\mu}^*G_{\mu}c_{\mu}^{-1}\gamma,
\nonumber
\\
\chi &= \Im\zeta,
\nonumber
\\
\tilde{V}_{a,b}(c) &=  \chi -d\eta_{\Psi'}(\chi) + \tau_{\Psi'}(\chi).
\nonumber
\end{align}
In the following we will compute the error  if $\mathcal{F}(a,b)\neq 0$.
We will use the formulas from Lemma \ref{lemma-ACIDs}.
To make the notation less cumbersome, in the following calculations we will write $\dbar$  for $\dbar_{[\mu]}$ and $N$ for $N_{a,b}$.

We first note that $dJ_{a,b}d$ still annihilates the exact term $d\eta_{\Psi'}(\chi)$, 
and by Lemma \ref{lemma-6smallerror} $dJ_{a,b}d\tau_{\Psi'}(\chi)$ is a small error term. Thus we have 
\begin{align*}
dJ_{a,b}d \tilde{V}_{a,b}(c)
=
dJ_{a,b}d \chi + 
\text{small error}.
\end{align*}
The mixed terms $(\partial+\dbar)Ji_N$ and $i_N J(\partial+\dbar)$ only contribute small error terms.
Therefore, we have
\begin{align*}
dJ_{a,b}d\chi
&=
i
(\partial\dbar\chi -\dbar\partial\chi)
+
\text{small error}
\\
&=
2i \partial\dbar\chi
+
\text{small error}
\\
&=
-2i \dbar\partial \chi 
+
\text{small error}.
\end{align*}
Thus we have
\begin{align}
dJ_{a,b}d\chi
=
\partial\dbar\zeta 
+
\dbar\partial\bar{\zeta}
+
\text{small error}.
\label{error-chi}
\end{align}
Under our assumption $\mathcal{H}^{1,2}=0$ the action of $[\dbar_{\mu},G_{\mu}]$ on $\Omega^{1,2}$ by Lemma\ref{lemma-commutator} only contributes small error terms. Therefore, we have
\begin{align}
\dbar\zeta
=
\gamma - \dbar^* c_{\mu}G_{\mu}c^{-1}_{\mu}\dbar\gamma + \text{small error}.
\label{error1-zeta}
\end{align}
Now we have
\begin{align}
\gamma=\beta_{1,2}-\partial\lambda-i_{N''}\bar{\lambda}.
\label{gamma}
\end{align}
$\rho=d\beta$ is equivalent to
\begin{subequations}
\begin{gather}
\dbar\beta_{3,0}+\partial\beta_{2,1}+ i_{N'}\beta_{1,2}=0,
\label{rho1}
\\
\dbar\beta_{2,1}+\partial\beta_{1,2}+i_{N''}\beta_{3,0}+i_{N'}\beta_{0,3}=\rho,
\label{rho2}
\\
\partial\beta_{0,3}+\dbar\beta_{1,2}+i_{N''}\beta_{2,1}=0.
\label{rho3}
\end{gather}
\end{subequations}
Therefore, from \eqref{gamma} and \eqref{rho3} we get
\begin{align*}
\dbar\gamma
=
\dbar\beta_{1,2}-\dbar\partial\lambda
-
\dbar i_{N''}\bar{\lambda}
=
-\partial\beta_{0,3} + \partial\dbar\lambda
+
\text{small error}.
\end{align*}
Just as in the integrable case we still have
$\dbar\lambda = \beta_{0,3}$ 
and thus 
\begin{align}
\dbar\gamma
=
\text{small error}.
\label{dbar-gamma}
\end{align}
Combining \eqref{error1-zeta} and \eqref{dbar-gamma} gives
\begin{align}
\dbar\zeta
=
\gamma
+ \text{small error}.
\label{dbar-zeta2}
\end{align}
By \eqref{gamma} and the identities from Lemma \ref{lemma-ACIDs} we get
\begin{align*}
\partial\gamma
=
\partial\beta_{1,2}-\partial^2\lambda
-
\partial i_{N''}\bar{\lambda}
=
\partial\beta_{1,2}+\text{small error},
\quad
\dbar\bar{\gamma}
=
\dbar\beta_{2,1}+\text{small error},
\end{align*}
and therefore by \eqref{rho2}
\begin{align}
\partial\gamma+\dbar\bar{\gamma}=\rho+\text{small error}.
\label{gamma-rho}
\end{align}
From \eqref{dbar-zeta2}
\begin{align*}
\partial\dbar\zeta
=
\partial\gamma
+\text{small error},
\\
\dbar\partial\bar{\zeta}
=
\dbar\bar{\gamma}+\text{small error},
\end{align*}
and thus with \eqref{gamma-rho}
\begin{align}
\partial\dbar\zeta 
+
\dbar\partial\bar{\zeta}
=
\rho
+
\text{small error}.
\label{zeta-rho}
\end{align}
Combining \eqref{error-chi} and \eqref{zeta-rho} gives
\begin{align*}
dJ_{a,b}d\chi
=
\rho
+
\text{small error}.
\end{align*}
Thus we have
\begin{align*}
\tilde{L}_{a,b}\tilde{V}_{a,b}(c)
=&
c
+
\text{small error}.
\end{align*}
This completes the proof of Theorem \ref{Thm-specialCase}.

\subsection{Proof under general hypothesis}

In this subsection we complete the proof for the general case of Theorem \ref{MainThm} where $\mathcal{K}_{\Psi_0}=0$ but $\mathcal{H}^{1,2}\oplus\mathcal{H}_M\neq 0$. We use the notation 
\begin{align*}
L = L_{0,0},
\quad
V=V_{0,0}.
\end{align*}
As we explained in the introduction, our boundary value problem has a variational
formulation. This implies that the linearised operator $L_0$ at a solution should be formally
self-adjoint; that is, we expect that the condition $\mathcal{K}_{\Psi_0} = 0$, which is that the kernel of $L_0$
is zero, should be equivalent to the vanishing of the cokernel. To make this precise we
need:

\begin{lemma}
\label{lemma-image}
The image of $L$ is closed and of finite codimension.
\end{lemma}
\begin{proof}
Let $c\in C$ and solve as in \eqref{2step-process}.
However, if $\mathcal{H}^{1,2}\neq 0$, then with Proposition \ref{Prop-UniSolveH0} and Lemma \ref{lemma-commutator} 
we now get
\begin{align*}\gamma = \dbar\zeta + H_0(\gamma).
\end{align*}
Then we have
\begin{align*}
c = *dJd\chi + f^Z(H_0(\gamma)),
\end{align*}
where $f^Z:\mathcal{H}^{1,2}\rightarrow C$ is the map $f^Z(x) = *\Re (\partial x)$.

The process of adapting $\chi$ to get a $2$-form which vanishes on $M$ similarly defines $y\in\mathcal{H}_M$ and a map $f^M:\mathcal{H}_M\rightarrow C$ such that
\begin{align*}
*dJdV_0 c = c + f^Z(H_0(\gamma)) + f^M(y).
\end{align*}

The map $f:= f^Z + f^M$ maps $\mathcal{H}^{1,2}\oplus \mathcal{H}$ to a finite dimensional subspace $E$ of $C$.
By the above discussion there is a continuous map $h: C\rightarrow \mathcal{H}^{1,2}\oplus \mathcal{H}$ such that 
\begin{align*}
L V c = c + f(h(c))
\end{align*}
for all $c\in C$. Thus, we have $ \im(\mathbbm{1}+f\circ h) \subset \im\, L$. The statement follows because $\im (\mathbbm{1}+f\circ h)$ is closed and of finite codimension.
\end{proof}

\begin{lemma}
\label{lemma-selfAdj}
$L_{a,b}:A\rightarrow C$ is self-adjoint with respect to the $L^2$-inner product for all $(a,b)\in U$.
\end{lemma}
\begin{proof}
We first note the identity
\begin{align}
\int_Z dJ_{a,b}d\sigma\wedge\zeta = \int_Z \sigma\wedge dJ_{a,b}d\zeta
\label{integrationID}
\end{align}
for $2$-forms  $\sigma$ and $\zeta$ which restrict to $0$ on the boundary.
This follows from Stokes' theorem and a simple calculation which shows that $J_{a,b}\gamma\wedge \mu = - \gamma\wedge J_{a,b}\mu$ for all $\gamma,\mu\in\Omega^3(Z)$. Let $\alpha,\beta\in A$. Substituting $\sigma = \alpha+\Gamma\alpha$ and $\zeta=\beta+\Gamma\beta$ 
in \eqref{integrationID} gives 
\begin{align*}
(*dJ_{a,b}d(\alpha+\Gamma\alpha),\beta+\Gamma\beta) = (\alpha+\Gamma\alpha,*dJ_{a,b}d(\beta+\Gamma\beta)).
\end{align*}
The operator $R$ from \eqref{def-phi} satisfies $\mathbbm{1}-R = d^*dG_6\pi_6$. The statement of the Lemma follows because if $\alpha\in A$ and $\tau\in\Omega^2_6$ with $\tau\bdry=0$ we have
\begin{align*}
(d^*d\tau,\alpha+\Gamma\alpha)
=
(\tau, \pi_6 d^*d(\alpha+\Gamma\alpha))=0,
\end{align*}
as $d(\alpha+\Gamma\alpha)\in T$.
\end{proof}

\begin{corollary}
$\mathcal{K}_{\Psi_0} = \ker L$ is a finite dimensional space. 
If $\mathcal{K}_{\Psi_0}=0$, then $L$ is surjective.
\end{corollary}

Choose a complement $C'$ to $E\subset C$ so $C= C'\oplus E$ and we have projections $\pi_{C'}:C\rightarrow C'$ and $\pi_E: C\rightarrow E$. Let $Y = L^{-1}(E)\subset A$. For each $(a,b)\in U$ we have
a linear map $\lambda_{a,b}:Y\rightarrow E$ which is tame in $(a,b)$ defined by $\lambda_{a,b}(y) = \pi_EL_{a,b}y$. By definition $\lambda_{0,0}$ is an isomorphism. We have
\begin{align*}
\lambda_{0,0}^{-1}\lambda_{a,b} - \mathbbm{1} 
=
\lambda_{0,0}^{-1}\pi_E(L_{a,b}-L_{0,0}),
\end{align*} 
which is $(a,b)$-tamely small of order $0$ by Lemmas \ref{lemma-comp-tamely-small} and \ref{lemma-adaptedMoserEstimate} and thus by Proposition \ref{prop-NeumannSeries} after shrinking $U$ we can assume that $\lambda_{a,b}$ has a tame inverse.

Similarly as in the proof of \ref{lemma-image},
let $c\in C$, and solve as in \eqref{2step-process}.
By Proposition \ref{Prop-UniSolveH0}
we now have
\begin{align*}
\gamma = c_{\mu}(\dbar_{\mu}\dbar^*_{\mu}G_{\mu} + \dbar^*_{\mu}\dbar_{\mu}G_{\mu} + H_{\mu})c_{\mu}^{-1}\gamma
=
\dbar_{[\mu]}\zeta
+
c_{\mu}\dbar^*_{\mu}[\dbar_{\mu},G_{\mu}]c_{\mu}^{-1}\gamma + \text{small error}
+
c_{\mu}H_{\mu}c_{\mu}^{-1}\gamma
.
\end{align*}
In the following write $\mathcal{O}^{\mathrm{ts}}_{a,b}(\sigma)$ for any $(a,b)$-tamely small operator (as defined in Definition \ref{def-mu-small}) applied to $\sigma$.
By Lemma \ref{lemma-commutator} the commutator $[\dbar_{\mu},G_{\mu}]$ contributes
a small error term and a term which is of the form $\mathcal{O}^{\mathrm{ts}}_{a,b}(H_{\mu}c_{\mu}^{-1}\gamma)$.
Thus we have
\begin{align*}
\gamma
=
\dbar_{\mu}\zeta
+
c_{\mu}H_{\mu}c_{\mu}^{-1}\gamma
+
\mathcal{O}^{\mathrm{ts}}_{a,b}( H_{\mu}c_{\mu}^{-1}\gamma)
+ \text{small error}
.
\end{align*}
Now there is a map $f^Z_{a,b}$ such that
\begin{align*}
c 
= 
R*dJ_{\Psi'}d\chi + f^Z_{a,b'}H_{\mu}c_{\mu}^{-1}\gamma 
+
\mathcal{O}^{\mathrm{ts}}_{a,b}( H_{\mu}c_{\mu}^{-1}\gamma)
+ 
\text{small error}
.
\end{align*}
Continuing as in the proof of Lemma \ref{lemma-image} we get maps 
$f^M_{a,b}:\mathcal{H}_M\rightarrow C$, $f_{a,b}:= f^Z_{a,b}+f^M_{a,b}$ and $h_{a,b}:C\rightarrow \mathcal{H}^{1,2}\oplus \mathcal{H}_M$ such that we have
\begin{align*}
L_{a,b}V_{a,b}c
=
c + f_{a,b}(h_{a,b}(c))
+ 
\mathcal{O}^{\mathrm{ts}}_{a,b}(h_{a,b}(c))
+
\text{small error}.
\end{align*}
The map $f_{a,b}-f$ is $(a,b)$-tamely small. E.g. the component $f^Z_{a,b} - f^Z$ is the composition of a tame map with $c_{\mu}\partial_{\mu}-\partial$ and thus by Lemmas \ref{lemma-comp-tamely-small} and \ref{lemma-adaptedMoserEstimate} is $(a,b)$-tamely small.
Therefore, we have 
\begin{align*}
L_{a,b}V_{a,b}c
=
c + f(h_{a,b}(c))
+ 
\mathcal{O}^{\mathrm{ts}}_{a,b}(h_{a,b}(c))
+
\text{small error}.
\end{align*}
Let $S_{a,b}=\lambda_{a,b}^{-1}\circ f\circ h_{a,b}:C\rightarrow Y$ and 
$\tilde{L}_{a,b}=\pi_{C'}\circ L_{a,b}:A\rightarrow C'$.
We have 
\begin{align*}
L_{a,b}(V_{a,b}c-S_{a,b}c) = c - \tilde{L}_{a,b}S_{a,b}c 
+ 
\mathcal{O}^{\mathrm{ts}}_{a,b}(h_{a,b}(c))
+
\text{small error}.
\end{align*}
By definition $\tilde{L}_{0,0}=0$, and thus we have
\begin{align*}
\tilde{L}_{a,b} = \tilde{L}_{a,b} - \tilde{L}_{0,0} = \pi_{C'}\circ (L_{a,b}-L_{0,0}).
\end{align*}
By Lemma \ref{lemma-adaptedMoserEstimate}  this is $(a,b)$-tamely small. 
To sum up, there is a $(a,b)$-tamely small map $T_{a,b}$ such that
\begin{align*}
L_{a,b}(V_{a,b}c-S_{a,b}c) = c - T_{a,b}(h_{a,b}(c))
+
\text{small error}.
\end{align*}
By Lemma \ref{lemma-comp-tamely-small} (iii) the term $T_{a,b}(h_{a,b}(c))$ is $(a,b)$-tamely small of order $0$, and thus 
by Proposition \ref{prop-NeumannSeries} the map $\mathbbm{1}-T_{a,b}\circ h_{a,b}$ has a tame inverse. 
We obtain a tame map
\begin{align*}
\tilde{V}_{a,b} = (V_{a,b}-S_{a,b})\circ (\mathbbm{1}-T_{a,b}\circ h_{a,b})^{-1}
\end{align*}
with
\begin{align}
L_{a,b}\tilde{V}_{a,b}c = c + \text{small error}.
\label{final-equ}
\end{align}
Therefore, we can apply Theorem \ref{Nash-Moser-Zehnder-Thm}
to the family of maps $\phi_b = \mathcal{F}(\,\cdot\,,b)$ to obtain the existence statement in Theorem \ref{MainThm}.

Finally we prove the uniqueness statement in Theorem \ref{MainThm}. If $\mathcal{F}(a,b)=0$, then by \eqref{final-equ}
$L_{a,b}$ is surjective, and thus by Lemma \ref{lemma-selfAdj}
injective. The statement follows from Theorem \ref{Thm-uniqueness}.

\section{The deformation space}

\label{section-defspace}

\subsection{The space of infinitesimal deformations}

Recall that $\cK$ is the space of infinitesimal deformations of a solution in a fixed enhanced boundary class:
$$  \cK=\frac{\{ d\sigma: \sigma\in \Omega^{2}(Z), dJd \sigma=0, \sigma\vert_{M}=0\}}{ \{ \mathcal{L}_{w}\Psi: w\vert_{M}=0\}} $$
Introduce a larger space
$$  \hcK=\frac{\{ \rho\in \Omega^{3}(Z): d\rho= 0, dJ\rho=0, \rho\vert_{M}=0\}}{\{ \mathcal{L}_{w}\Psi: w\vert_{M}=0\}}.  $$
This is the space of infinitesimal deformations ignoring the enhancement. So $\cK$ is the kernel of a map $\Pi_{Z}: \hcK\rightarrow H^{3}(Z,M)$. In terms of the decomposition
$\rho= \rho_{3,0} + \rho_{2,1} + \overline{\rho_{2,1}} + \overline {\rho_{3,0}}$ the equations $d\rho=dJ\rho=0$ are equivalent to 
\begin{align}
\db \rho_{2,1}=0,\quad \partial \rho_{2,1}= \db \rho_{3,0}.
\label{inf-def-eq}
\end{align}

Recall some linear algebra at a point $p$  in our submanifold $M\subset Z$. Let $z:TZ_{p}\rightarrow \bC$ be a complex linear map such that $M_{p}$ is the kernel of
${\rm Re}(z)$. For each $(p,q)$ we have a space of \lq\lq normal'' forms $\Lambda^{p,q}_{N}\subset \Lambda^{p,q}_{Z}=\Lambda^{p,q}T^{*}Z_{p}$. These are the forms $d\overline{z}\wedge \mu_{p,q-1}$. The quotient space is 
$$  \Lambda^{p,q}_{T}= \Lambda^{p,q}_{Z}/\Lambda^{p,q}_{N}. $$
The restriction map $\Lambda^{p,q}_{Z}\rightarrow \Lambda^{p,q}(H)$ factors through $\Lambda^{p,q}_{T}$. Let $L^{p,q}\subset \Lambda^{p,q}_{T}$ be the kernel of $\Lambda^{p,q}_{T}\rightarrow \Lambda^{p,q}_{H}$ so we have an exact sequence
\begin{align}
0\rightarrow L^{p,q}_{T}\rightarrow \Lambda^{p,q}_{T}\rightarrow \Lambda^{p,q}(H) \rightarrow 0.
\label{seq-forms-bdry}
\end{align}
The space $L^{p,q}_{T}$ consists of forms in $\Lambda^{p,q}_{Z}$ which can be written as $dz\wedge \mu_{p-1,q}$ modulo those which can be written as
$dz d\overline{z}\wedge \mu_{p-1,q-1}$. Since $dz d\overline{z}$ vanishes on $TM$ there is a well-defined map from $L^{p,q}_{T}$ to $\Lambda^{p+q}T^{*}M$. 
A  choice of contact form, and so Reeb vector field, defines a splitting of \eqref{seq-forms-bdry} and we can then write
  $$  \Lambda^{p,q}_{T}= L^{p,q}_{T} \oplus \Lambda^{p,q}(H)= \Lambda^{p-1,q}(H)\wedge dz \oplus \Lambda^{p,q}(H). $$
but unlike the subspace $L^{p,q}_{T}$ this splitting is not canonical. 

The $\dbar$-complex on $Z$ induces the $\dbar_{b}$-complex on sections  $\Omega^{*,*}_{T}$ of $\Lambda^{p,q}_{T}$.  Given a choice of splitting, the $\dbar_{b}$ operator is the sum of $\dbar_{H}: \Omega^{p,q}_{H}\rightarrow \Omega^{p,q+1}_{H}$ and algebraic operators
\begin{align*}
\omega: \Omega^{p,q}_{H}\wedge dz \rightarrow \Omega^{p+1,q+1}_{H},\quad\quad S:\Omega^{p,q}_{H}\rightarrow \Omega^{p-1,q+1}_{H}\wedge dz.
\end{align*}
Here we are writing $dz$ for the sum $I\theta+i\theta$ where $\theta$ is the contact form on $M$.

With this background in place we return to our specific situation. By definition $\cH$ is the set of real $\gamma\in  \Omega^{1,1}_{H}$ with $d_{H}\gamma=0, \omega\wedge \gamma=0$. These conditions are equivalent to $d(\gamma\wedge \theta)=0$.  For $\gamma\in \cH$ the section  
$\gamma\wedge dz$ of $ L^{2,1}_{T}$ is $\db_{b}$-closed, so we have a map 
\begin{align}
\iota: \cH\rightarrow H^{2,1}_{T}.
\label{def-iota}
\end{align}

\begin{lemma}
\label{lemma-3-form-bdry}
Suppose $\rho\in \Lambda^{3}(T^{*}Z)$ and $\rho\bdry=0$. Then the image of $\rho_{2,1}$ in $\Lambda^{2,1}_{T}$ lies in $L^{2,1}_{T}$ and the restriction of $J\rho$ to $M$ is $2i$ times the image of this under the map $L^{2,1}_{T}\rightarrow \Lambda^{3}T^{*}M$. 
\end{lemma}
\begin{proof}
The condition that $\rho\vert_{TM}=0$ means that we can write
$$  \rho= (dz+d\overline{z})(\sigma_{2,0}+ \sigma_{1,1}+ \sigma_{0,2}). $$
Then $\rho_{2,1}= dz\sigma_{1,1}+ d\overline{z} \sigma_{2,0}$
which equals $dz\sigma_{1,1}$ modulo $\Lambda^{2,1}_{N}$, so $\rho_{2,1}$ maps to $L^{2,1}_{T}$. Applying the definition of $J$:
$$  J\rho= idz\sigma_{2,0} + i( dz \sigma_{1,1}+ d\overline{z} \sigma_{2,0}) - i( dz\sigma_{0,2}+ d\overline{z}\sigma_{1,1}) - id\overline{z}\sigma_{0,2}, $$
which is $i(dz+d\overline{z})(\sigma_{2,0}-\sigma_{0,2}) + i(dz-d\overline{z}) \sigma_{1,1}$ and the result follows.
\end{proof}

Now let $\rho$ be a representative of a class in $\hcK$ so $\rho$ is a real $3$-form with $d\rho=dJ\rho=0$ and $\rho\vert_{M}=0$. The restriction of $J\rho$ to $M$ is closed and by the Lemma, lies in 
$\cH\wedge \theta$. This gives a map 
\begin{align*}
b:\hcK\rightarrow \cH.
\end{align*}
If $[\rho]\in\ker\, b$, then $\rho_{2,1}\in\Omega^{2,1}_N$,
and by \eqref{inf-def-eq} we get a map
\begin{align*}
\lambda:       {\rm ker}\, b\rightarrow  H^{2,1}_{N}.
\end{align*}  
The composition of $\iota$ with the coboundary map of the long exact cohomology sequence gives a map
\begin{align*}
\mu: \cH\rightarrow H^{2,2}_{N}.
\end{align*}
Because $\iota \circ b: \hcK \rightarrow H^{2,1}_{T}$
maps $\rho$ to the image of $\rho_{2,1}$ in $H^{2,1}_{T}$,
we have $\mu\circ b = 0$.

\begin{proposition}
\label{prop-exact-seq}
Suppose that $H^{3,1}_{Z}=0$ and that any function $f$ on $M$ with $\db f=0$ is the restriction of a holomorphic function on $Z$. Then $\lambda$ is an isomorphism and the image of $b$ is equal to ${\rm ker}\ \mu$. 
\end{proposition}

So under these hypotheses we have an exact sequence
    $$    H^{2,1}_{N}\rightarrow \hcK \rightarrow \cH \rightarrow H^{2,2}_{N}. $$
If the hypotheses are not satisfied the same arguments as below  will give  more complicated descriptions of the kernel and image of $b$ for example if $H^{3,1}\neq 0$, but with the same condition on holomorphic functions the kernel of $b$ is isomorphic to the kernel of
$$ \partial: H^{2,1}_{N}\rightarrow H^{3,1}. $$

\begin{proof}[Proof of Proposition \ref{prop-exact-seq}]
{\bf $\lambda$ is injective:} Suppose $[\rho]\in \hcK$ with $b(\rho)=0$ and $\rho_{2,1}=0$ in $H^{2,1}_{N}$. Thus $\rho_{2,1}=\db \eta_{2,0}$ with $\eta_{2,0}\in \Omega^{2,0}_{N}$. This means that $\eta_{2,0}= i_{w} \Psi$ for a vector field $w$ vanishing on the boundary and  $\rho- d (\eta_{2,0}+\overline{\eta_{2,0}})$ defines the same class in $\hcK$. Thus we can suppose that $\rho_{2,1}=0$. So $\rho = \rho_{3,0}+\overline{ \rho_{3,0}} $ for a holomorphic $3$-form $\rho_{3,0}= f \Psi$. Now the condition that  the restriction of $\rho$ to $M$ is zero implies that the holomorphic function $f$ vanishes on the boundary, so is zero everywhere.

\

{\bf $\lambda$ is surjective:}
Now suppose that $\sigma_{2,1}\in \Omega^{2,1}_{N}$ with $\db \sigma_{2,1}=0$. Since $H^{3,1}=0$, by assumption, we can find a $\sigma_{3,0}$ such that $\partial \sigma_{2,1}= - \db \sigma_{3,0}$. The fact that $\sigma_{2,1}\in \Omega^{2,1}_{N}$ implies that the restriction of  $\sigma= \sigma_{2,1}+ \sigma_{3,0}$to the boundary has the form  $f \theta \wedge (\alpha+ i\beta)$ for a complex-valued function $f$. Now $\sigma$ is closed on $Z$ so 
on $M$ we have
$$   d (f \theta(\alpha+i\beta)=0$$
which is equivalent to $\db_{b}f=0$. By assumption, $f$ extends to a holomorphic function on $Z$ and changing $\sigma_{3,0}$ to $\sigma_{3,0} - f \Psi$ we can suppose that the resriction of $\sigma$ to the boundary is $0$. Then $\sigma+ \overline{\sigma}$ gives an element of ${\rm ker} b$ mapping to the class of $\sigma_{2,1}$ in $H^{2,1}_{N}$.

{\bf ${\rm Im}\  b = {\rm ker}\ \mu$:}
Suppose that $\gamma \in \cH$ and $\mu(\gamma)=0$. The long exact cohomology sequence implies that we can find a $\sigma_{2,1}\in \Omega^{2,1}$mapping to $\gamma\wedge dz\in \Omega^{2,1}_{T}$ with $\db\sigma_{2,1}=0$. As above, the hypothesis $H^{3,1}=0$ means that we can find  $\sigma_{2,0}$ so that
$d\sigma=0$ where $\sigma= \sigma_{2,1} + \sigma_{3,0}$. Now he restriction of $\sigma$ to the boundary has the form $ i \gamma \wedge \theta + f (\alpha+i\beta)\wedge \theta$  and the fact that $d_{H} \gamma=0$ implies again that $\db_{b}f=0$, from which the argument is the same as before. 
\end{proof}

When $Z$ is  a Stein manifold the hypotheses in Proposition \ref{prop-exact-seq} are satisfied and the cohomology groups $H^{2,1}_{N}$ and $H^{2,2}_{N}$ vanish,  so we simply get
\begin{corollary}
If $Z$ is a Stein manifold, then $\hcK \cong \Pi_{Z}:\cH\rightarrow H^{3}(Z,M)$.
\end{corollary}

There is a related discussion for the embedding problem considered in the previous paper \cite{embedding-paper}. There we defined a map $\Pi_M: \mathcal{H}\rightarrow H^3(M,\bR)$ taking $h\in\mathcal{H}$ to the class of $\theta\wedge h$ and we showed that the obstruction space for the embedding problem is the kernel of $\Pi_M$.

\subsection{Comparison with CR deformation theory and extension to higher dimensions}

There is  a large literature on deformations of CR structures and the problems of realising such deformations by change of embedding in a fixed ambient complex manifold or and deformations of complex manifolds with boundary. The authors do not  have enough knowledge of this literature to give a reasonable survey but we just make some remarks. 

The infinitesimal deformations of a {\it closed} complex manifold $Z$ are given by the sheaf cohomology group $H^{1}(TZ)$. In the case of a Calabi-Yau 3-fold we have an isomorphism $TZ=\Lambda^{2}T^{*}Z$ so $H^{1}(TZ)= H^{2,1}$. Similarly the discussion of deformations of a general CR structure on a hypersurface $M$ involve   cohomology groups with co-efficients in the bundle $TZ\vert_{M}$ to $M$ but in our situation everythiong can be expressed in terms of differential forms.  In the tangential complex the elements of $\Omega^{2,1}_{T}$ correspond to deformations of the \lq\lq almost CR'' structure, {\it i.e.} a field of subspaces $H\subset TM$ and complex structure on $H$. In the representation
$$\Omega^{2,1}_{T}= \Omega^{2,1}_{H} \oplus \Omega^{1,1}_{H}\wedge \theta, $$
the first summand corresponds to deformations of the subspace $H$ and the second to deformations of complex structure on $H$. The kernel of $\db_{b}:\Omega^{2,1}_{T}\rightarrow \Omega^{2,2}_{T}$ gives solutions of the linearisation of the integrability conditions for the almost CR structure. In the representation
$$  \Omega^{2,0}_{T}= \Omega^{2,0}_{H}\oplus \Omega^{1,0}_{H}\wedge \theta, $$
the second term corresponds to vector fields on $M$ which take values in $H\subset TM$. The first term can be identified with complex valued functions on $M$. The real-valued functions give vector fields on $M$, multiples of the Reeb field $v$ and the pure imaginary valued functions give normal vector fields, multiples of $Iv$. So $\Omega^{2,0}_{T}$ can be viewed as sections of $TZ\vert_{M}$, which are infinitesimal deformations of the embedding of $M$ in $ Z$. The operator
$\db_{b}:\Omega^{2,0}_{T}\rightarrow \Omega^{2,1}_{T}$ gives the infinitesimal action of these deformations on the CR structures,  so the cohomology $H^{2,1}_{T}$ represents deformations of the CR structure modulo those deformations which can be realised by deforming $M$ in $Z$.  (To study deformations of the CR structure modulo diffeomorphisms of $M$ one would have to restrict to elements of $\Omega^{2,}_{T}$ with real component in $\Omega^{2,0}_{H}$.) 

The exact sequence (under our two hypotheses)
    $$   H^{2,1}_{N}\rightarrow \hcK \rightarrow \cH \rightarrow H^{2,2}_{N} $$
maps term-by-term to the long exact sequence in $\db$-Neumann theory
$$   H^{2,1}_{N}\rightarrow H^{2,1}_{Z}\rightarrow H^{2,1}_{T}\rightarrow H^{2,2}_{N} . $$
   The map $\iota:\cH\rightarrow H^{2,1}_{T}$ gives the infinitesimal deformation of the CR structure defined by the deformation of the complex $3$-form $\Phi$ on $M$.  
The second sequence encodes, at the formal infinitesimal level, the deformation problem for complex manifolds with boundary. A class in $H^{2,1}_{T}$ lifts to $H^{2,1}_{Z}$ if the infinitesimal deformation of the CR structure can be realised by an infinitesimal deformation of $Z$.  For genuine deformations, Kiremidjian \cite{Kiremidjian1979ADE} showed that if
$H^{2,2}_{N}=0$ then any deformation of the CR structure is obtained from one of $Z$ (this is the formulation of Kiremdijian's result in the Calabi-Yau case).
For  a Stein manifold $Z$,  all the cohomology groups in the second sequence vanish but in the first sequence $\cH$ and $\hcK$ can be non-trivial.  

In part,  Hitchin's work \cite{Hitchin} gave a new point of view on the known deformation theory of closed Calabi-Yau threefolds. Similarly, our work in this paper gives, in part, a new point of view on known deformation results for CR structures. If one knows that any deformation of the $3$-form $\phi=\alpha\wedge\theta$ on $M$ has a matching deformation of $\tilde{\phi}=\beta\wedge \theta$---which gives a deformation of the CR structure---then the existence of the deformation of $Z$, under suitable cohomology vanishing conditions, could be obtained from the result of Kiremidjian. Other relevent references, include \cite{Akahori2002DeformationTO} (which is focused on the 5-dimensional case)
and many papers of Kuranishi, Akahori, Miyajima and other authors.

Our definition of the space $\cH$ extends to $CR$ manifolds in all dimensions. As recalled above, a CR structure defines subspaces $L^{p,q}\subset \Lambda^{p,q}_{T}$ which can also be regarded as $(p+q)$-forms on $M$. 
We can define
$$  \cH_{\bC}^{p,q}= \{ \chi\in L^{p,q}_{T}\subset \Omega^{p+q}(M): d\chi = 0\}. $$
When $M$ has dimension $5$ and $p=2, q=1$ this is the complexification of our space $\cH$.
In general, we have maps $\cH_{\bC}^{p,q}\rightarrow H^{p+q}(M,\bC)$ and $\cH_{\bC}^{p,q}\rightarrow H^{p,q}_{T}$. If the Levi form is nondegenerate, these spaces are zero for rather trivial reasons for small $p,q$.
The interesting case seems to be for a pseudonvex structure  on  a manifold $M$ of dimension $4k+1$ with $p=k+1, q=k$. Fixing a choice of contact form we have
$$ \cH_{\bC}^{k+1,k}\equiv \{\gamma\in \Omega_{H}^{k,k}: \omega \wedge \gamma=0, d_{H}\gamma=0\}. $$
The algebraic condition $\omega \wedge \gamma=0$ defines the primitive subspace $P^{k,k}_{H}\subset \Omega_{H}^{k,k}$ which consists of anti-self-dual forms, so $d{H}\gamma=0$ implies $d^{*}_{H}\gamma=0$. Arguments similar to those in \cite{embedding-paper} show that $\cH_{\bC}^{k+1,k}$ is finite-dimensional.  We do not at present know how this fits into the  literature on CR geometry, or what the geometric significance of these spaces might be.

\subsection{An example: $T^{3}\times B^{3}$}

Take $\bC^{3}= \bR^{3}+ i \bR^{3}$ in the standard way and let $W$ be the product of the unit ball $B^{3}$ in $\bR^{3}$ with $i\bR^{3}$. Let $Z$ be the quotient of
$W$ by the lattice $i\bZ^{3}\subset i\bR^{3}$. So $Z$ is a domain in $\bC^{3}/i\bZ^{3}$ with boundary $M=S^{2}\times T^{3}$.  We take the holomorphic $3$-form $\psi+ i\tpsi= i dz_{1} dz_{2} dz_{3}$ so
$$\psi= dy_{1}dy_{2}dy_{3}- \sum_{\rcyc} dy_{i} dx_{j}dx_{k} \  \  \  \  \  \tpsi= dx_{1}dx_{2}dx_{3}- \sum_{\rcyc} dx_{i}dy_{j}dy_{k}. $$

We begin by finding the space $\cH_{M}$. The contact form is $\theta= \sum x_{i} dy_{i}$, so $\omega=d\theta= \sum dx_{i} dy_{i}$.  We have $\psi\vert_{M}= \alpha \wedge \theta, \tpsi\vert_{M}= \beta\wedge \theta$ where
$$\alpha= \sum_{\rcyc} x_{i} (dy_{j} dy_{k}- dx_{j}dx_{k})  \ \ \ \ \ \ \ \ , \ \ \ \ \  \beta=\sum_{\rcyc} x_{i} (dx_{j}dy_{k}- dx_{k} dy_{j}) . $$
The self-dual forms on the contact $4$-planes are spanned by $\alpha,\beta, \omega$. The rotations $SO(3)$ act on the cover $\partial W= S^{2}\times i\bR^{3}$ and the forms $\alpha, \beta,\omega$ are $SO(3)$-invariant. There is another $SO(3)$- invariant form
$$\gamma= \sum_{\rcyc} x_{i} (dx_{j}dx_{k}+ dy_{j} dy_{k}). $$
This is an anti-self dual form. Set
$$  \lambda= dy_{1}dy_{2} dy_{3} +\sum_{\rcyc} dy_{i} dx_{j} dx_{k}. $$
We have
$$\gamma\wedge \theta = \sum_{a}\sum_{\rcyc} x_{a} x_{i} dy_{a}(dx_{j}dx_{k}+ dy_{j}dy_{k})$$  and a moments thought shows that this is equal to the restriction of $\lambda$ to the 2-sphere.  So $d(\gamma\wedge \theta)=0$ and $\gamma$ is in our space $\cH_{M}$.
\begin{proposition}
The space $\cH_{M}$ is one-dimensional, spanned by $\gamma$. 
\end{proposition}

The contact subspace $H\subset TM$ splits as $TS^{2} \oplus I TS^{2}$ 
We have $\Lambda^{-}_{H}= \bR \gamma \oplus Q$ where $Q$ is the orthogonal complement of $\gamma$.  This can be identified with the bundle $s^{2}_{0}$ trace-free quadratic forms on $TS^{2}$, lifted to $M$. (The identification goes by mapping $q\in Q$ to the quadratic form $w\mapsto q(w, I w)$.) 
 Recall that on $M$ we have our complex
   \begin{equation}\Omega^{0}_{H}\stackrel{d_{H}}{\rightarrow} \Omega^{1}_{H} \stackrel{d^{-}_{H}}{\rightarrow} \Omega^{-}_{H}. \end{equation}
We consider first the forms that are invariant under translations by $i\bR^{3}$. This becomes a complex on $S^{2}$ and a short calculation identifies this as the sum of the de Rham complex
$$   \Omega^{0} \rightarrow \Omega^{1}_{X} \rightarrow \Omega_{S^{2}}^{2}, $$
and $$    0\rightarrow \Omega^{1}_{Y}\stackrel{D}{\rightarrow}  \Gamma(s^{2}_{0}). $$
Here we are writing $\Omega^{1}_{X}, \Omega^{1}_{Y}$ for the sections of the duals of $TS^{2}, ITS^{2}$ respectively (so they give two copies of the $1$-forms on $S^{2}$) and we are identifying $2$-forms on $S^{2}$ with 
$\bR\gamma\subset \Lambda^{-}_{H}$ in the obvious way, using the standard area form on $S^{2}$. The operator $D$ maps a $1$-form $\eta$ to the symmetric, trace-free, component of the covariant derivative. This can also be identified with $\db$-operator on vector fields.

The operator $D$ is surjective so we see immediately that the cokernel of $d^{-}_{H}$ on these translation-invarient forms is the $1$-dimensional space spanned by $\gamma$. To handle the general case we use a Fourier decomposition in the $T^{3}$ factor. We complexify the vector bundles and for $\xi$ in $ i \bR^{3}$ we consider forms of the shape
$a(x) {\rm exp}( \xi.y)$. For each $\xi$ we get a complex over $S^{2}$ which differs from that above (when $\xi=0$) by three extra terms. The vector $\xi$ defines a $1$-form $\tilde{\xi}= \sum \xi_{i} dx_{i}$ on $S^{2}$. The extra terms are:
\begin{itemize} \item $A:\Omega^{0}\rightarrow \Omega^{1}_{Y}$ with $A(f)= f \tilde{\xi}$;
\item $B:\Omega^{1}_{X}\rightarrow \Gamma(s^{2}_{0})$ with $B(\eta)= -\tilde{\xi}*\eta$, where $*$ is the composite of the product map $T^{*}\times T^{*}\mapsto s^{2}$ and projection to the trace-free part. 
\item  $C:\Omega^{1}_{Y}\rightarrow \Omega^{2}$ with $C(\eta)=  \tilde{\xi}\wedge \eta$. 
\end{itemize}
Proposition * amounts to the statement that the cohomology in degree $2$ of this complex vanishes, for all non-zero $\xi$. In general, associated to our complex * we have a Laplace operator $\Delta_{H}^{-}= d_{H}^{-} \left( d_{H}^{-}\right)^{*}$ on $\Omega^{-}_{H}$.  In our situation, for a fixed Fourier component $\xi$, we get an operator over $S^{2}$
$$  \Delta_{\xi}: \Omega^{2}\oplus \Gamma(Q)\rightarrow \Omega^{2}+ \Gamma(Q). $$
We claim that this is diagonal with respect to the direct sum decomposition, that is
$$   D C^{*}+ B d^{*} = 0 \ \ \ \ \ \   C D^{*}+ d B^{*}=0 . $$
The second equation is the adjoint of the first, so it suffices to prove the first.  Write $\mu$ for the standard area form on $S^{2}$. Then 
$$   C^{*}(f \mu) = f j(\tilde{\xi}), $$
where $j$ is the complex structure on $T^{*}S^{2}$,  regarded as a complex line  bundle in the usual way.  (The signs are confusing here, recall that $\xi$ is pure imaginary and this gives a change in sign in the adjoint.)
Now $D(j(\tilde{\xi}))=0$ so
$DC^{*} (f \mu)= df* (j(\tilde{\xi}))$. On the other hand $d^{*}(f\mu)= j(df)$ so $B d^{*} =-\tilde{\xi}*j(df)$. Now the claim follows from the identity, for any $1$-forms $\eta_{1}, \eta_{2}$
$$  \eta_{1}* j(\eta_{2})= \eta_{2}*j(\eta_{1}) . $$

Given this claim it is clear that $\Delta_{\xi}$ is a strictly positive operator, for non-zero $\xi$. 

Now we have
\begin{proposition}
The maps $\Pi_{M}: \cH_{M}\rightarrow H^{3}(M)$ and $\Pi_{Z}: \cH_{M}\rightarrow H^{3}(Z,M)$ are injective.
\end{proposition}

For the first part  we just have to show that the closed form $\gamma\wedge \theta$ is non-zero in $H^{3}(M)$  but this is immediate since $\gamma\wedge \theta= dy_{1}dy_{2}dy_{3} - \sum_{a}\sum_{\rcyc} x_{a} x_{i} dy_{a} d_{j} dx_{k}$ which has non-zero integral over ${\rm pt.}\times T^{3}\subset S^{2}\times T^{3}$.

For the second part we need to find a lift of $\gamma$ to $\hcK$. Set $\chi=dx_{1}dx_{2}dx_{3}$ so 
$$ 8  J(\chi)=  J((dz_{1}+ d\overline{z}_{1})(dz_{2}+d\overline{z}_{2})(dz_{3}+ d\overline{z}_{3})= {\rm Im}(dz_{1}dz_{2}dz_{3} +\sum_{\rcyc} d\overline{z}_{i}dz_{j}dz_{k}). $$
A few lines of calculation give
$$   J(\chi)= \frac{1}{2} ( dy_{1}dy_{2}dy_{3} + \sum_{\rcyc} dy_{i}dx_{j}dx_{k})= \frac{1}{2}\lambda. $$
Then $d\chi=dJ(\chi)=0$,  the restriction of $\chi$ to $S^{2}$ vanishes and the restriction of $2J\chi$ to the 2-sphere is $\gamma\wedge\theta$. So $2\chi$ defines a class in $\hcK$ with $b(2\chi)= \gamma$. The relative homology group $H_{3}(Z,M)$ is generated by $B^{3}\times {\rm pt}$ and clearly the integral of $\chi$ over this ball is non-zero.

Another point of view on the calculations above is given by the discussion of the \lq\lq constant mean curvature'' condition in \cite{papertwo} .

The consequence of  the second part of  Proposition * is that $Z$ is rigid so Theorem 1 applies. The consequence of the first part of the Proposition is that any small exact deformation of the boundary form can be realised by a small deformation of $M$ in $ \bC^{3}/i \bZ^{3}$.

\appendix

\section{Nash--Moser background}

\label{appendix-Nash-Moser}

In the following let $L(m)f$ be a linear operator, where the variables $m$ and $f$ are sections of some vector bundles either over a closed manifold or a compact manifold with boundary. Let $U$ be a neighborhood of $m=0$.

\begin{definition}
\label{def-tame}
We say that $L(m)$ is \textit{tame (of order s)} if there exists r such that for each $n$ we have an estimate
\begin{align*}
\|L(m)f\|_n
\lesssim
\|f\|_{n+s} + (1+\|m\|_{n+r}) \|f\|_s
\end{align*}
for all $m\in U$ and all $f$.
\end{definition}

\begin{definition}
\label{def-mu-small}
Suppose that there are $r$ and $s$ such that for each $n$ and $\epsilon>0$ there exists $C_n(\epsilon)$ such that $L(m)$ for each $m\in U$ and each $f$ satisfies an estimate
\begin{align*}
\|L(m)f\|_n
\leq
\epsilon
\|m\|_r \|f\|_{n+s} + C_n(\epsilon) \|m\|_{n+r} \|f\|_s.
\end{align*}
Then we say that $L(m)$ is $m$\textit{-tamely small (of order $s$)}.
\end{definition}
A straightforward application of the definitions gives

\begin{lemma}
\label{lemma-comp-tamely-small}
Suppose that for $m\in U$ we have linear operators 
\begin{align*}
K(m):\Gamma(V_1)\rightarrow \Gamma(V_2), \quad
L(m):\Gamma(V_2)\rightarrow \Gamma(V_3).
\end{align*}
Then we have:
\begin{compactenum}[(i)]
\item If $K(m)$ is tame of order $r$ and $L(m)$ is $m$-tamely small of order $t$, then $L(m)\circ K(m)$ is $m$-tamely small of order $r+t$.
\item
If $K(m)$ is $m$-tamely small of order $r$ and $L(m)$ is tame of order $t$, then $L(m)\circ K(m)$ is $m$-tamely small of order $r+t$.
\item In cases (i) and (ii), if for all $m\in U$ $K(m)$ maps to a finite dimensional subspace $H\subset \Gamma(V_2)$, then $L(m)\circ K(m)$ is $m$-tamely small of order $0$.
\end{compactenum}
\end{lemma}
\begin{proof}
(i) and (ii) are a straightforward consequence of the definitions.
To prove (iii), if for example we are in case (i), then because on the finite dimensional vector space $H$ all norms are equivalent, $K(m)$ actually is tame in any order, in particular tame of order $-t$.
\end{proof}

An important example of an $m$-tamely small operator is a differential operator $L(m)$ which vanishes for $m=0$.

\begin{lemma}[Improved Moser estimate 4]
\label{lemma-adaptedMoserEstimate}
Suppose that $L(m)f$ is a partial differential operator, possibly nonlinear of degree $r$ in $m$ and linear of degree $s$ in $f$. Suppose $L(0)f=0$ for all $f$. Then for each $n\in\bZ$ and $\epsilon>0$ there is a constant $C_n(\epsilon)$ such that
\begin{align*}
\|L(m)f\|_n
\leq 
\epsilon
[[m]]_r \|f\|_{n+s}
+
C_n(\epsilon)
[[m]]_{n+r}
\|f\|_s.
\end{align*}
\end{lemma}
\begin{proof}
We have \cite[p. 440, ``Moser estimate 4'']{Hamilton1977Deformation2}
\begin{align*}
\|L(m)f\|_n
\lesssim
\sum_{i+j=n}
([[m]]_r \|f\|_{n+s})^{\frac{j}{n}}
([[m]]_{n+r}\|f\|_s)^{\frac{j}{n}}.
\end{align*}
The statement follows if  we apply the inequality
\begin{align*}
ab \leq (sc) a^{\frac{1}{p}}+(lc)b^{\frac{1}{q}},
\end{align*}
where $a$ and $b$ are positive numbers and $p+q=1$. We apply this to each term, setting $a=[[m]]_r \|f\|_{n+s}$, $b=[[m]]_{n+r}\|f\|_s$,
$p=\frac{j}{n}$ and $q=\frac{i}{n}$.
\end{proof}

For perturbations of the identity by an $m$-tamely small operator of order $0$ we can construct a Neumann series to get a tame inverse.
\begin{proposition}
\label{prop-NeumannSeries}
Suppose there exists $r$ such that for each $n$ there exists $C_n$ such that
\begin{align}
\|L(m)f\|_n \leq \|m\|_r \|f\|_n + C_n \|m\|_{n+r} \|f\|_r
\label{goodEstimate}
\end{align}
for all $m\in U$ and all $f$ (this is in particular true if $L(m)$ is $m$-tamely small of order $0$). Then the operator $P(m)f=f-L(m)f$ has an inverse $P^{-1}(m)f$ which satisfy a tame estimate
\begin{align}
\|P^{-1}(m)f\|_n 
\lesssim \|f\|_n + \|m\|_{n+r} \|f\|_r.
\label{tame-inverse}
\end{align}
\end{proposition}
\begin{proof}
Using \eqref{goodEstimate} with $n=r$ gives
\begin{align}
\|L(m)f\|_r \leq C \|m\|_{2r} \|f\|_r
\label{n=0}
\end{align}
with $C=1+C_r$.
Next we will show that for all $k$ we have
\begin{align}
\|L^k(m)f\|_n
\leq 
\|m\|_r^k \|f\|_n
+ C_n \|m\|_{n+r} \|f\|_r
\sum_{j=0}^{k-1}
\|m\|_r^{k-1-j} (C \|m\|_{2r})^j.
\label{iteration-estimate}
\end{align}
The case $k=1$ is the inequality \eqref{goodEstimate}.
Assume that \eqref{iteration-estimate} is true for $k-1$. Then 
\eqref{goodEstimate} and \eqref{n=0} give
\begin{gather*}
\|L^k(m)f\|_n \leq \|m\|_r \|L^{k-1}(m)f\|_n + C_n \|m\|_{n+r} \|L^{k-1}(m)f\|_r
\\
\leq
\|m\|_r \|L^{k-1}(m)f\|_n + C_n \|m\|_{n+r} \|f\|_r (C\|m\|_{2r})^{k-1}.
\end{gather*}
\eqref{iteration-estimate} follows by induction.
If $\|m\|_{2r} < \min(1/4,1/(4C))$, then
\begin{align*}
\|L^k(m)f\|_n 
\leq 4^{-k}
\|f\|_n
+ 4 C_n \|m\|_{n+r} \|f\|_r\, k 4^{-k}.
\end{align*}
Because $2^{-k}k\rightarrow 0$ as $k\rightarrow \infty$, there exists $K_n$ such that for all $k\geq K_n$ we have
\begin{align}
\|L^k(m)f\|_n 
\leq 4^{-k}
\|f\|_n
+ 2^{-k} \|m\|_{n+r} \|f\|_r.
\label{sumableIneq}
\end{align}
By comparison with the geometric series, the series $\sum_{k=0}^{\infty} L^k(m)f$ converges in every $L^2_n$ to a smooth limit.
Taking the limit of \eqref{sumableIneq} proves \eqref{tame-inverse}. 
\end{proof}

\begin{theorem}
\label{Nash-Moser-Zehnder-Thm}
Let $\phi:\mathcal{C}^{\infty}(Z,V)\rightarrow \mathcal{C}^{\infty}(Z,W)$ be an operator between sections of Hermitian vector bundles $V$ and $W$ over a compact manifold $Z$ (possibly with boundary).
Denote the $L^2_s$-norms for sections of $V$ and $W$ by $|\cdot|_s$
 and $\|\cdot\|_s$, respectively. Suppose that there exist $u_0\in\mathcal{C}^{\infty}(Z,V)$, an integer $d>0$, a real number $\delta$ and constants $C_1, C_2$ and $(C_s)_{s\geq d}$ such that for any $u,v,w\in\mathcal{C}^{\infty}(Z,V)$,
 \begin{align*}
 |u-u_0|<\delta
 \Rightarrow
 \begin{cases}
 \forall s\geq d, \|\phi(u)\|_s \leq C_s (1+|u|_{s+d}),
 \\
 \|\phi'(u)v\|_{2d} \leq C_1 |v|_{3d},
 \\
 \|\phi''(u)(v,w)\|_{2d}
 \leq 
 C_2 |v|_{3d}|w|_{3d},
 \end{cases}
 \end{align*}
 where $\phi'$ and $\phi''$ are the first and second derivative of $\phi$, respectively.
Moreover, suppose that for every $u\in\mathcal{C}^{\infty}(Z,V)$
which satisfies $|u-u_0|_{3d}<\delta$ there exist operators $\psi(u):\mathcal{C}^{\infty}(Z,W)\rightarrow \mathcal{C}^{\infty}(Z,V)$ and $\mathcal{E}(u):\mathcal{C}^{\infty}(Z,W)\times \mathcal{C}^{\infty}(Z,W)\rightarrow \mathcal{C}^{\infty}(Z,W)$ and a constant $C_3$ such that for any $\varphi,\rho\in\mathcal{C}^{\infty}(Z,W)$ we have
\begin{gather*}
\phi'(u)\psi(u)\varphi = \varphi + \mathcal{E}(u)(\phi(u),\varphi),
\\
\|\mathcal{E}(u)(\varphi,\rho)\|_{2d}
\leq C_3 \|\varphi\|_{3d}\|\rho\|_{3d},
\\
\forall s\geq d,
|\psi(u)\varphi|_{s}
\leq 
C_s(\|\varphi\|_{s+d}+|u|_{s+d}\|\varphi\|_{2d}).
\end{gather*}  
Then if $\|\phi(u_0)\|_{2d}$ is sufficiently small (with respect to some upper bound on $1/\delta$ and a finite number of the constants $C_s$), there exists $u\in\mathcal{C}^{\infty}(Z,V)$ such that $\phi(u)=0$.
\end{theorem}
This is an adaption of the presentation of the Nash--Moser implicit function theorem in \cite{StRaymond}. We briefly describe how the proof must be adapted.
The solution $u$ is found as the limit of a certain sequence $(u_k)_{k\geq 0}$. To estimate $\|\phi(u_k)\|_{2d}$, the Taylor expansion of order 1 is used. The calculation on p. 222 in \cite{StRaymond} now gives
\begin{align*}
\phi(u_{k+1})
=
\varphi_1 
+ \varphi_2
+\varphi_3,
\end{align*}
where the additional term
\begin{align*}
\varphi_3
=
-\mathcal{E}(u_k)(\phi(u_k),\phi(u_k))
\end{align*}
 incorporates the failure of $\psi(u)$ to be a  right inverse of $\phi'(u)$.
We need to bound the $L^2_{2d}$-norm of this term by $C_0 \theta_k^{-6}$.
Using the estimate for the error term gives
\begin{align*}
\|\mathcal{E}(u_k)(\phi(u_k),\phi(u_k))\|_{2d}
\leq C_3 \|\phi(u_k)\|_{3d}^2.
\end{align*}
Now we do not have a bound on $\|\phi(u_k)\|_{3d}$, but the bound $\|\phi(u_k)\|_{2d}\leq \theta_k^{-4}$. 
Setting $\bar{\theta}_k = \theta_k^{1/d}$, for $t\geq 3d$ we can use the smoothing operator $S_{\bar{\theta}_k}$ to interpolate
\begin{gather*}
\|\phi(u_k)\|_{3d}
\leq 
\|S_{\bar{\theta}_k}\phi(u_k)\|_{3d}
+
\|(\mathbbm{1}-S_{\bar{\theta}_k})\phi(u_k)\|_{3d}
\\
\leq 
C_{3d,2d} \bar{\theta}_k^d \|\phi(u_k)\|_{2d}
+
C_{3d,t} \bar{\theta}_k^{3d-t} \|\phi(u_k)\|_t
\\
\leq 
C_{3d,2d} \bar{\theta}_k^d \theta_k^{-4}
+
C_{3d,t} C_t \bar{\theta}_k^{3d-t}  (1+|u_k|_{t+d})
\\
=
C_{3d,2d} \theta_k^{-3}
+
C_{3d,t} C_t \bar{\theta}_k^{3d-t} (1+|u_k|_{t+d})
.
\end{gather*}
Now $t+d > 3d$, so we cannot use the bound $|u_k-u_0|_{3d} < \delta$.
However, as Saint Raymond explains \cite[p. 223, proof of Lemma 2,]{StRaymond} the sequence $(1+|u_k|_s)\theta_k^{-N}$ with $N=4(2d+1)$ is monotone decreasing for a fixed $s\geq 2d$, so the second term can be estimated as
\begin{gather*}
C_{3d,t} C_t \bar{\theta}_k^{3d-t} (1+|u_k|_{t+d})
=
C_{3d,t} C_t \bar{\theta}_k^{3d-t} \theta_k^N (1+|u_k|_{t+d})\theta_k^{-N}
\\
\leq 
C_{3d,t} C_t \bar{\theta}_k^{3d-t} \theta_k^N (1+|u_0|_{t+d})\theta_0^{-N}.
\end{gather*}
If we set $t=3d+d(N+3)$, then this is bounded by $C_0 \theta_k^{-3}$, which is what we needed.

\begin{theorem}\cite[Theorem 3.3.3]{HamiltonIFT}
\label{Thm-uniqueness}
Let $A, B, C$ be tame Frechet spaces and
\begin{align*}
\mathcal{F}: U \subset (A\times B) \rightarrow C
\end{align*}
a smooth tame map with $\mathcal{F}(0,0)=0$. 
Suppose that there is a smooth tame map
\begin{align*}
V: U \times C \rightarrow A
\end{align*}
which is a right inverse for $D_{\sigma}\mathcal{F}(\sigma,\xi)$ whenever $\mathcal{F}(\sigma,\xi)=0$, and that $D_{\sigma}\mathcal{F}(\sigma,\xi)$ is injective whenever $\mathcal{F}(\sigma,\xi)=0$. 
Then there is a neighborhood $U'$  of $0$ in $A$ and a neighborhood $U''$ of $0$ in $B$, such that if $\sigma_1$ and $\sigma_2$ are in $U'$ and $\xi$ is in $U''$ with $\mathcal{F}(\sigma_1,\xi)=0$ and $\mathcal{F}(\sigma_2,\xi)=0$, then $\sigma_1=\sigma_2$. 
\end{theorem}
\begin{proof}
The hypothesis imply that $V(\sigma,\xi)$ is also a left inverse for
$D_{\sigma}\mathcal{F}(\sigma,\xi)$ whenever $\mathcal{F}(\sigma,\xi)=0$, so that we can write
\begin{align*}
V(\sigma,\xi)D_{\sigma}\mathcal{F}(\sigma,\xi)\zeta
=
\zeta + Q(\sigma,\xi)\{\mathcal{F}(\sigma,\xi),\zeta)
\end{align*} 
with $Q(\sigma,\xi)\{\mathcal{F}(\sigma,\xi),\zeta)=0$ for all $\zeta$ if $\mathcal{F}(\sigma,\xi)=0$. The statement follows as in the proof of \cite[Theorem 3.3.3]{HamiltonIFT}.
\end{proof}

\bibliographystyle{amsplain}

\bibliography{bibliografie}

\end{document}